\documentclass{amsart}
\usepackage{amsmath}%
\usepackage{amsfonts}%
\usepackage{amssymb}%
\usepackage{graphicx}
\usepackage{amsthm}
\usepackage{euscript}
\usepackage{latexsym}
\usepackage[all]{xy}
\usepackage{verbatim}
\newtheorem{thm}{Theorem}[section]
\newtheorem{defn}[thm]{Definition}
\newtheorem{prop}[thm]{Proposition}

\newtheorem{cor}[thm]{Corollary}
\newtheorem{lemma}[thm]{Lemma}

\newcommand{\field}[1]  {\mathbb{#1}} % Use blackboard bold for these sets

\newcommand{\N}         {\field{N}}
\newcommand{\Z}         {\field{Z}}

\newcommand{\Q}         {\field{Q}}
\newcommand{\F}         {\field{F}}

\begin{document}
\title{On the 2-typical de Rham-Witt complex}
\author{Viorel Costeanu}
\date{July 16, 2007}
\subjclass{Primary 13K05; Secondary 19D55} %
\keywords{de Rham-Witt, topological cyclic homology, algebraic K-theory}%

\begin{abstract}
In this paper we introduce the 2-typical de Rham-Witt complex for arbitrary commutative, unital rings and log-rings. We describe this complex for the rings $\Z$ and $\Z_{(2)}$, for the log-ring $(\Z_{(2)},M)$ with the canonical log-structure, and we describe its behaviour under polynomial extensions. In an appendix we also describe the $p$-typical de Rham-Witt complex of $(\Z_{(p)},M)$ for $p$ odd.\end{abstract}
\maketitle

\section{Introduction}

The $p$-typical de Rham-Witt complex was introduced by Bloch, Deligne, and Illusie for $\F_p$-algebras (see \cite{B}, \cite{Ill}). The definition was generalized by Hesselholt and Madsen to $\Z_{(p)}$-algebras, for $p$ odd (see \cite{HM1}, \cite{HM2}, \cite{He1}). Hesselholt and Madsen's construction was motivated by the effort to understand $\operatorname{TR},$ an  object that appears in algebraic topology and is related to topological cyclic homology and to higher algebraic $K$-theory. More precisely, for a fixed prime $p$ and a $\Z_{(p)}$-algebra $A$, one defines:
\[ \operatorname{TR}^n_q(A;p) = \pi_q(T(A)^{C_{p^{n-1}}}),\]
where $T(A)$ is the topological Hochschild spectrum associated to $A,$ and $C_r\subset S^1$ is the cyclic group of order $r.$ As $n$ and $q$ vary these groups are related by certain operators $F, V, R, d, \iota$ which satisfy several relations. One notes that $\iota$ is induced by the multiplication with the element $\eta\in\pi_1^{s}{S^0}$ from stable homotopy. This element has order 2, so the operator $\iota$ is trivial if 2 is invertible. This is the case if $A$ is a $\Z_{(p)}$-algebra with $p$ odd, and this explains why the case $p=2$ is different from $p$ odd. 

A first step in understanding $\operatorname{TR}$ is to understand the universal example of an object that has the same algebraic structure as $\operatorname{TR}.$ The algebraic structure of $\operatorname{TR}$ is captured by the notion of a Witt complex, that we will give shortly. The fact that $\operatorname{TR}$ is a Witt complex was proved by Hesselholt in \cite{He1}. Before giving the definition we make precise what we mean by a pro-object and a strict map of pro-objects. We let $\Z$ be the category associated with the poset $(\Z, \geq);$ a pro-object in a category $\mathcal{C}$ is a covariant functor $X:\Z\to \mathcal{C},$ in other words a sequence of objects $\{X_n\}_{n\in\Z}$ and of morphisms $R:X_n\to X_{n-1}.$ A strict map of pro-objects is a natural transformation of functors, that is a sequence of maps $f_n:X_n\to Y_n$ that commutes with the maps $R.$

\begin{defn}A 2-typical Witt complex consists of:
\begin{itemize}
\item[(i)] a graded-commutative pro-graded ring $\{E_n^*,\;R: E_n^* \to E_{n-1}^*\}_{n\in \Z},$ such that $E_n^*=0$ for all $n\leq 0$. The index $n$ is called the level. 
\item[(ii)] a strict map of pro-rings 
$\lambda : W_{\bullet} (A) \to E_{\bullet} ^0$
from the pro-ring of Witt vectors of A.
\item[(iii)]a strict map of pro-graded rings 
$$ F : E^*_{\bullet} \to E^*_{{\bullet} -1},$$
such that $\lambda F = F \lambda.$ 
\item[(iv)]a strict map of pro-graded $E_{\bullet}^*$-modules
$$V:F_*E_{\bullet}^* \to E_{{\bullet}+1}^*$$
such that $\lambda V =V\lambda$ and $FV=2$. 
The linearity of $V$ means that $V(x)y=V(xF(y)), \forall x\in E_n^*, y\in E_{n+1}^*$.
\item[(vi)] a strict map of pro-graded abelian groups $d: E^*_{\bullet} \to 
E^{*+1}_{\bullet}$, which is a derivation, in the sense that 
$$d(xy)=d(x)y+(-1)^{\deg(x)}xd(y)$$

The  operator $\iota : E^*_{\bullet} \to 
E^{*+1}_{\bullet}$  is by definition multiplication by the element $\frac 
{d\lambda[-1]_n}{\lambda[-1]_n}$, where  $[a]_n =(a, 0, \dots, 0) \in W_n(A)$ is the multiplicative representative. 

The operators $F,$ $V,$ $d,$ and $\iota$ are required to satisfy the following relations: 

$$FdV = d+\iota,$$
$$dd=d\iota=\iota d,$$
$$ Fd\lambda([a]_n) = \lambda([a]_{n-1})d\lambda([a]_{n-1}]), \text{ for all } a\in A.$$
\end{itemize}
\end{defn}

Visually, a Witt complex is a two dimensional array: 

\[\xymatrix@R=9ex@C=5em{&&&&&\\
    &W_2(A)\ar@{.}[u]\ar[d]^R \ar[r]^\lambda\ar@<-1ex>[d]_F
    &E^0_2\ar[d]^R\ar@{.}[u]\ar@<-.5ex>[r]_\iota\ar@<.5ex>[r]^d\ar@<-1ex>[d]_F
    &E^1_2\ar[d]^R\ar@{.}[u]\ar@<-.5ex>[r]_\iota\ar@<.5ex>[r]^d\ar@<-1ex>[d]_F
    &E^2_2\ar[d]^R\ar@{.}[u]\ar@{.}[r]\ar@<-1ex>[d]_F&\\
          A\ar[uur]^{[-]_n}\ar[ur]\ar[r]^{[-]_1}
          &W_1(A)\ar[r]^\lambda\ar@<-1ex>@/_/[u]_V
          &E^0_1\ar@<-.5ex>[r]_\iota\ar@<.5ex>[r]^d\ar@<-1ex>@/_/[u]_V
          &E^1_1\ar@<-.5ex>[r]_\iota\ar@<.5ex>[r]^d\ar@<-1ex>@/_/[u]_V
          &E^2_1\ar@<-1ex>@/_/[u]_V\ar@{.}[r]&\\}
\]

A map of 2-typical Witt complexes is a map $f:E_{\bullet}^* \to {E^\prime}_{\bullet}^*$ 
of pro-graded rings such that $\lambda^\prime = f\lambda,$ $fd=df,$
$F^\prime f=fF,$ and $V^\prime f=fV.$

The paper is organized as follows. In Section 2 we discuss Witt vectors, the de Rham complex, and Witt complexes in general. We also derive the identity that expresses the Teichm\"{u}ller representative of an integer as a combination of a system of generators:
\[[a]_n=a[1]_n+\sum_{i=1}^{n-1}\frac{a^{2^i}-a^{2^{i-1}}}{2^i}V^i[1]_{n-i}.\]

In the third section we prove, using category theory, that the category of 2-typical Witt complexes over a given ring admits an initial object, and that by definition is the de Rham-Witt complex of the ring. A similar result holds for the more general notion of a log-ring. 

Section 4 contains several calculations. The first result of this section is the structure theorem of the de Rham-Witt complex of the ring of rational integers $\Z.$ It states that in degree zero it is the pro-ring of Witt vectors of the integers, in degree one it is generated by the elements $dV^i(1),$ and in degrees above one it vanishes: 
\begin{eqnarray}
W_n\Omega^0_\field{Z}&\cong& \bigoplus_{i=0}^{n-1}\field{Z}\cdot V^i(1),\\
W_n\Omega^1_\field{Z}&\cong& \bigoplus_{i=1}^{n-1}\field{Z}/{2^i\field{Z}}\cdot dV^i(1) ,\\
W_n\Omega^i_\field{Z}&=&0,\quad \mbox{for } i\geq 2.
\end{eqnarray}
The product rule and the action of the operators are given in Theorem \ref{structureforZ} below. We note that, additively, the formula for  the 2-typical de Rham-Witt complex is similar to the one for the  $p$-typical de Rham-Witt complex for $p$ odd. The differences appear in the product rule and the action of the operators $d,$ $F,$ and, of course, $\iota.$  In a remark at the end of the section we note that a very similar result holds for the de Rham-Witt complex of the ring $\Z_{(2)}.$

In Section 4 we also describe the behaviour of the de Rham-Witt complex under polynomial extensions. Again the result is similar to the one in the $p$-typical case, for $p$ odd, which is found in Section 4.2 in \cite{HM1}. The de Rham-Witt complex of the ring $A[X]$ consists of formal sums of four types of elements:
\begin{itemize}
\item Type 1: elements of the form $a[X]^j$, where $a\in W_n\Omega_A^q$,
\item Type 2: elements of the form $b[X]^{k-1}d[X]$, where $b\in W_n\Omega_A^{q-1}$,
\item Type 3: elements of the form $V^r(c[X]^l)$, where $r>0$, $c\in W_{n-r}\Omega_A^q$, and $l$ is odd,
\item Type 4: elements of the form $dV^s(e[X]^m)$, where $s>0$, $e\in W_{n-s}\Omega_A^{q-1}$, and $m$ is odd.
 \end{itemize}
The product rule and the action of the operators are given explicitely. 

In the last part of the fourth section we define the notion of a Witt complex for log-rings and we compute the 2-typical de Rham-Witt complex of the log-ring $(\Z_{(2)},M),$ where $M=\Q^*\cap\Z_{(2)}\hookrightarrow \Z_{(2)}$ is the canonical log-structure. The difference from the 2-typical de Rham-Witt complex of $\Z_{(2)}$ and of $(\Z_{(2)},M)$ is the element $d\operatorname{log}[2]:$
\begin{align}
W_n\Omega^0_{(\Z_{(2)},M)}&\cong W_n\Omega^0_{\Z_{(2)}}\cong W_n({\Z_{(2)}}),\\
 W_n\Omega^1_{(\Z_{(2)},M)}&\cong W_n\Omega^1_{\Z_{(2)}}\oplus \Z/{2^n\Z}d\operatorname{log}[2]_n,\\
W_n\Omega^i_{(\Z_{(2)},M)}&=0, \text{for all } i\geq 2.
\end{align}
An interesting formula in this context is: 
\[V(d\operatorname{log}[2]_n)= 2d\operatorname{log}[2]_{n+1}+dV[1]_{n}-dV^2[1]_{n-1}+4dV^3[1]_{n-2}.\]

The paper has two appendices. In the first one we describe the structure of the $p$-typical de Rham-Witt complex of the log ring $(\Z_{(p)},M),$ with $p$ odd. This result is very similar to the one for $p=2$, the difference being in the product formulas and the action of the operator $V.$ We note here that there are two distinct cases for $p$ odd, namely $p=3$ and $p\geq 5.$ For example, the mentioned formula becomes:
\[V(d\operatorname{log}[p]_n)=
\begin{cases}
3d\operatorname{log}[3]_{n+1}+dV[1]_n+3dV[1]_{n-1},& \text{if } p=3,\\
p\,d\operatorname{log}[p]_{n+1}+dV[1]_n,& \text{if } p\geq 5.
\end{cases}\]

In the second appendix, which is rather technical, we verify the associativity of the multiplication defined in Section 4, subsection 4.2.   

In this paper all rings are associative, commutative, and unital. Graded rings are graded commutative, or anti-symmetric, meaning that, for every two elements $x,y$ of degrees $|x|, |y|,$ respectively, one has
\[xy=(-1)^{|x||y|}yx.\]

{\bf Acknowledgement.} This paper is based on the author's Ph.D. dissertation written under the direction of Lars Hesselholt at MIT. The author wants to thank Lars Hesselholt for his enthusiasm and inspiring guidance.

\section{Generalities: Witt vectors, the de Rham complex, and Witt complexes}

In this section we recall the Witt vectors and the de Rham complex. The 
standard references for these are \cite{Mu}, \cite{Mat}, respectively. 
Then we derive some elementary results for Witt vectors and Witt complexes.   

The de Rham complex of a ring $A$ is  the exterior algebra on the module 
of Kaehler differentials over $A$. More precisely, if $I$ is the kernel of the 
multiplication $A\otimes A \to A$, the module of Kaehler differentials 
is defined to be $ \Omega_A^1 = I/ {I^2}$; the map $d: A \to \Omega^1_A$ 
defined by $da = a\otimes 1- 1\otimes a +I^2$ is the universal derivation 
from $A$ to an $A$-module. 
The de Rham-complex is the exterior algebra 
\begin{align*}
\Omega^*_A &=\Lambda^*_A\Omega^1_A\\
\intertext{with differential}
d(a_0da_1\dots da_n)&=da_0da_1\dots da_n,\\
\intertext{where the exterior algebra of an $A$-module $M$ is} 
\Lambda^*(M) &=T_A(M)/\langle m\otimes m \mid m\in M \rangle.
\end{align*}

In this paper we will need a related construction, that of a universal anti-symmetric differential graded algebra over the ring $A$. By this we mean a graded algebra over $A$ which is commutative in the graded sense and is endowed with a $\Z$-linear differential of degree $1$, which is also a derivation. We will denote this by $\tilde{\Omega}^*_A$. Explicitly,
\[\tilde{\Omega}^*_A =\tilde{\Lambda}^*_A\Omega^1_A,\]
where:
\[\tilde{\Lambda}^*(M)=T(M)/\langle m\otimes n+n\otimes m \mid m, n\in M \rangle\]
is the universal anti-symmetric graded $A$-algebra generated by the $A$-module $M$.  When $2$ is invertible in $A$ the two constructions give the same result as the ideals $\langle m\otimes m \mid m\in M \rangle$ and $\langle m\otimes n+n\otimes m \mid m, n\in M\rangle$
are the same. In this paper we cannot assume that $2$ is invertible and this is why we need the second construction. 

The ring $W_n(A)$ of Witt vectors of length $n$ in $A$ is the set of 
$n$-tuples in $A$ with the following ring structure. One defines the 
``ghost'' map 
$$w:W_n(A) \to A^n$$
with components:
$$w_i(a_0,\dots,a_n)= a_0^{2^i}+ 2a_1^{2^{i-1}}+\dots 2^i a_i.$$
To add or multiply two vectors $a$ and $b$ one maps them via $w$ in 
$A^n$, adds or multiplies them componentwise, then uses $w^{-1}$ to map 
them back in $W_n(A)$. Of course one has to check that the sum or product 
of $w(a)$ and $w(b)$ are in  the image of the ghost map and that their
preimage is unique. That they are in the image follows from a lemma of Dwork; 
the uniqueness of the preimage is true only when $A$ has no $2$-torsion, 
which will be the case for the rings considered in this paper. When $A$ 
has $2$-torsion one has to give a canonical element in the preimage, and 
this is done requiring that the ghost map be a natural transformation of 
functors from rings to rings. 

The projection onto the first $n-1$ factors is a ring homomorphism
$$R: W_n(A) \to W_{n-1}(A),$$
called {\it restriction}, and this makes $W_{\bullet}(A)$ a pro-ring. There 
is a second ring homomorphism, the {\it Frobenius},
$$F:W_n(A)\to W_{n-1}(A),$$
$$w(F(a_0,\dots,a_{n-1}))=(w_1(a),\dots,w_{n-1}(a)),$$
and a $W_n(A)$-linear map, {\it Verschiebung},
$$V:F_*W_{n-1}(A) \to W_n(A)$$ 
$$V(a_0,\dots,a_{n-2})=(0,a_0,\dots,a_{n-2})$$
The notation $F_*W_{n-1}(A)$ indicates that $W_{n-1}(A)$ is considered a $W_n(A)$-module 
via the Frobenius map $F:W_n(A)\to W_{n-1}(A)$. The linearity of $V$ means 
therefore that $xV(y)=V(F(x)y),$ for all $x \in W_n(A)$ and $y \in W_{n-1}(A),$ 
formula known as Frobenius reciprocity. 
Both Frobenius and Verschiebung commute with the restriction maps. 
The {\it Teichm\"{u}ller} map is the multiplicative map
$$[\phantom{a}]_n: A \to W_n(A),$$
$$[a]_n=(a,0,\dots, 0).$$

We list now a few numerical results, some of which are not available in the odd 
prime case. 

\begin{prop}\label{relation-1}
In the ring of 2-typical Witt vectors of length n, $W_n(A)$,
$$[-1]_n= -[1]_n +V[1]_{n-1}.$$
\end{prop}
\begin{proof}
 It is enough to prove this relation for $A=\Z.$ In ghost coordinates, 
\begin{align*}w([-1]_n) &= (-1,1,\dots,1),\\
w([1]_n) &= (1,1,\dots,1),\\
w(V[1]_{n-1}) &= (0,2,\dots,2).
\end{align*}
The relation follows from the fact that addition is done component-wise in 
these coordinates and the ghost map is injective for $A=\Z.$ 
\end{proof}

\begin{prop}
In the ring of $2$-typical Witt vectors of length n, $W_n(\Z)$, there are $4$ 
square roots of unity, $[1]_n, [-1]_n, -[1]_n,-[-1]_n$. 
\end{prop}
\begin{proof}
Let $a = (a_0,\dots,a_{n-1}) \in W_n(A)$ be a square root of unity. Let 
$(w_0,\dots,w_{n-1})$
be its ghost coordinates. Then $(w_0^2, \dots, w_{n-1}^2) = (1,\dots, 1)$.
From $w_0^2=1$ we get $a_0^2=1$, hence $a_0 =\pm 1$. Equating the second ghost
coordinate we obtain: $ (a_0^2+2a_1)^2=1 \Rightarrow (1+2a_1)^2=1 \Rightarrow
a_1=0$ or $a_1=-1$. We will prove by induction that for $s\geq 1$, $a_s=a_1$. 
Assume this is true for $s-1$. We have two cases, $a_1=0$ and $a_1=-1$.
\begin{itemize}
\item[(i)]$a_1=0$: $w_s^2=(a_0^{2^s}+\cdots+2^{s-1}a_{s_1}^2 + 2^sa_s)^2=1$, 
so $(1+2^sa_s)^2=1$ and the unique integral solution is $a_s = 0$. 
 \item[(ii)]$a_1=-1$: $w_s^2=(a_0^{2^s}+\cdots+2^{s-1}a_{s_1}^2 + 2^sa_s)^2=1
\Rightarrow(1+2+\cdots+2^{s-1}+2^s a_s)^2=1 \Rightarrow (2^s-1+2^s a_s)^2=1$ 
and the unique integral solution is $a_s = -1$. 
\end{itemize}
Therefore the solutions of the equation $a^2=1$ are the vectors $(\pm 1, 0, 
\dots, 0)$ and $(\pm 1, -1,\cdots, -1)$. An examination of these vectors shows 
that they are exactly those listed in the statement. 
\end{proof}

\begin{prop}\label{wittbasis}
In the ring of $p$-typical Witt vectors of length $n$, $W_n(\field{Z})$, the vectors
$\{[1]_n, V([1]_{n-1}, \dots, V^{n-1}([1]_1)\}$ form a $\field{Z}$-basis. A vector 
$a=(a_0, \dots, a_{n-1})\in W_n(\Z)$ with ghost coordinates
$(w_0,\dots, w_{n-1})$ can be written in this basis as:
$$a=\sum_{s=0}^{n-1}c_sV^s([1]_{n-s}),$$
where 
$$c_s=\begin{cases} w_0 &\text{if $s=0$}\\
                    p^{-s}(w_s-w_{s-1})&\text{if $1\leq s \leq n-1$.}
      \end{cases} $$
The multiplication in this basis is given by the rule: 
$$V^i([1]_{n-i})V^j([1]_{n-j})=p^iV^j([1]_{n-j}), \quad \mbox{if}\quad i\leq j.$$
\end{prop}

\begin{proof}

In ghost coordinates, $V^s[1]_{n-s}=(0,\dots,0,p^s,\dots,p^s)$, the first 
$s$ coordinates being zero. Since the addition is component-wise it follows 
that these vectors are linearly independent. The multiplication is also component-wise
and the product formula follows.

We show that they form a system of generators. For a vector 
$a=(a_0, \dots, a_{n-1})\in W_n(\Z)$ with ghost coordinates
$(w_0,\dots, w_{n-1})$ we find the coefficients $c_i$ by induction. 
Equating the first ghost coordinate we get $c_0=w_0=a_0$. Assume we have found
$c_0,\dots,c_{s_1}$. We equate the $s$-th ghost coordinate 

$$w_s = \sum_{i=0}^s c_i p^i =\sum_{i=0}^{s-1} c_i p^i+c_s p^s =w_{s-1} + p^s 
c_s$$
and therefore, $c_s = p^{-s}(w_s-w_{s-1})$. These numbers are a priori rational. To finish the proof we need to show that they are integers. 

\begin{align*}
c_s&=p^{-s}(w_s-w_{s-1})\\
   &=p^{-s}(\sum_{i=0}^s p^ia_i^{p^{s-i}}-\sum_{i=0}^{s-1}p^ia_i^{p^{s-1-i}})\\
   &=p^{-s}(p^sa_s +\sum_{i=0}^{s-1} p^i(a_i^{p^{s-i}}-a_i^{p^{s-1-i}}))\\
   &=a_s+\sum_{i=0}^{s-1}p^{i-s}a_i^{p^{s-1-i}}(a_i^{p^{s-i}-p^{s-1-i}}-1)
\end{align*}
It remains to show that for every integer $a$ and every non-negative integer $n$:
\[a^{p^{n-1}}(a^{p^{n}-p^{n-1}}-1)\equiv 0\pmod{p^n} \]
There are two cases. If $v_p(a)\geq 1$, then $v_p(a^{p^{n-1}})\geq p^{n-1}\geq n$, and if $v_p(a)=0$, then $a^{p^{n}-p^{n-1}}-1=a^{\phi(p^n)}-1\equiv 0\pmod{p^n}.$
\end{proof}

\begin{cor}\label{formulafor2}
In the ring of $2$-typical Witt vectors of length n, $W_n(\Z)$, for every integer $a$, one has:
$$[a]_n=c_0[1]_{n}+c_1V[1]_{n-1} +\cdots+ c_{n-1} V^{n-1}[1]_{1},$$
where $c_0=a$ and $c_i=2^{-i}(a^{2^i}-a^{2^{i-1}})$.
\end{cor}

\begin{prop}
In every 2-typical Witt complex $E_{\bullet}^*$ the following relations hold:
\begin{align*}Vd &= 2dV, \\ 
              dF &= 2Fd, \\ 
           V(x)dV(y)&=V(xdy)+\iota V(xy) .
  \end{align*}
\end{prop}
\begin{proof}
We will use the relations from the definition of a Witt complex:
 $$\begin{array}{ccc}V(xF(y))=V(x)y, & FdV=d+\iota,& FV=2.
  \end{array}$$
We have:
\begin{eqnarray*} 
     Vd(x)&=&V(d+\iota+\iota)(x)=V(FdV+\iota)(x)\\
          &=&V(1)dV(x)+V\iota(x)\\
          &=& d(V(1)V(x))-dV(1)V(x)+V\iota(x)\\
          &=& d(V(FV(1)x))-V(FdV(1)x)+V\iota(x)\\
          &=& dV(2x)-V((d+\iota)(1)x)+V\iota(x)\\
          &=& 2dV(x)-V(d(1)x)-V(\iota x)+V(\iota x)\\
          &=& 2dV.\\
          &&\\
     dF(x)&=& (d+\iota)F(x)-\iota F(x)= FdVF(x) -\iota F(x)\\
          &=& Fd(V(1)x)-F\iota(x)=F(dV(1)x+V(1)d(x))-F\iota(x)\\
          &=& FdV(1)F(x)+FV(1Fdx)-F\iota(x)\\
          &=&(d+\iota)(1)F(x)+FVFdx-F\iota(x)\\
          &=& F(2dx)=2Fdx.\\
          &&\\
V(x)dV(y)&=& V(xFdV(y))\\
          &=& V(x(d+\iota)y)=V(xdy)+\iota V(xy).
 \end{eqnarray*}
\end{proof} 
\begin{prop}\label{iota}
$\iota ([1]_n)=\sum _{s=1}^{n-1} 2^{s-1}dV^s([1]_{n-s}).$
\end{prop}

\begin{proof} Since $2\iota([1]_n)=0$, we we'll prove that $\iota ([1]_n)=
-\sum _{s=1}^{n-1} 2^{s-1}dV^s([1]_{n-s}).$ 
The proof is by induction on $n$, starting with the case $n=1$ which is 
trivial. Assume the statement for $n-1$. We will use the relations  
$d([1]_n)=0$ and $([-1]_n)^2 = [1]_n$.
\begin{eqnarray*}
   [-1]_n &=& -[1]_n+V([1]_{n-1}),\\
d([-1]_n) &=& dV([1]_{n-1}),\\
\frac{d([-1]_n)}{[-1]_n} &=& \frac {dV([1]_{n-1})}{[-1]_n}=[-1]_n dV([1]_{n-1}),\\
\iota([1]_n) &=& [-1]_ndV([1]_{n-1})\\
             &=&(-[1]_n+V([1]_{n-1}))dV([1]_{n-1})\\
             &=&-dV([1]_{n-1})+V(FdV([1]_{n-1}))\\
             &=&-dV([1]_{n-1})+V((d+\iota)([1]_{n-1}))\\
             &=&-dV([1]_{n-1})+V(\iota([1]_{n-1}))\\
             &=&-dV([1]_{n-1})-\sum _{s=1}^{n-2} 2^{s-1}VdV^s([1]_{n-s}).
\end{eqnarray*}
The statement now follows from $Vd=2dV$. 
\end{proof}
 
\begin{prop}$\iota^2=0$.
\end{prop}

\begin{proof}
Again it is enough to prove $\iota^2([1]_n)=0$. We do this by induction. The 
case $n=1$ is trivial. Assume the statement for $n-1$. 
\begin{align*}
\iota^2 ([1]_n) &= \frac{d([-1]_n)}{[-1]_n}\frac{d([-1]_n)}{[-1]_n}\\
                &= (d([-1]_n))^2\\
                &= (d(-[1]_n+V([1]_{n-1}))^2\\
                &= dV([1]_{n-1})dV([1]_{n-1})\\
                &= d(V([1]_{n-1}) dV([1]_{n-1})) -V([1]_{n-1})ddV([1]_{n-1})\\
                &= d(V(FdV([1]_{n-1})))-V(FdV\iota([1]_{n-1}))\\
                &= d(V((d+\iota)([1]_{n-1})))-V((d+\iota)\iota([1]_{n-1})))\\
                &= dV\iota([1]_{n-1})) - V(\iota([1]_{n-1})\iota([1]_{n-1})).
 \end{align*}
The second summand is zero by induction. We show that the first summand is 
also zero. For this we will use the previous lemma and the relations $Vd=2dV,$ $2\iota =0,$ and $dd=d\iota$: 

\begin{align*} 
dV\iota([1]_{n-1})) &= dV(\sum _{s=1}^{n-2} 2^{s-1}dV^s([1]_{n-s})\\
                    &= \sum _{s=1}^{n-2} 2^sddV^{s+1}([1]_{n-s})\\
                    &= \sum _{s=1}^{n-2} (2^s\iota) d V^{s+1}([1]_{n-s})=0.   
\end{align*}
This completes the proof. 
\end{proof}

\section{The de Rham-Witt complex}
\subsection{Existence}
The Witt complexes over a ring $A$ form a category $\mathcal{W}_A$. 
We will 
prove that this category has an initial object. We call this object 
the de Rham-Witt complex of $A$ and  denote it    $W_{\bullet}\Omega^*_A$
To prove the existence of an initial 
object we use the Freyd adjoint functor theorem \cite[p.116]{Mac}. 

\begin{thm}\label{initialobject}
The category $\mathcal{W}_A$ of Witt complexes over $A$ has an initial 
object.
\end{thm}
\begin{proof} The category $\mathcal{W}_A$ has all small limits, so we need to prove that the  solution set condition is verified.
First we note that at each level a Witt complex is also a DG-ring. 
The differential is defined as follows: 

\begin{center}
 $D:E_{\bullet}^n\to E_{\bullet}^{n+1}, 
D= \begin{cases}
             d , &\text{if $n$ = even;} \\
             d+\iota ,& \text{if $n$ = odd.}
   \end{cases} $
\end{center}
\begin{prop}
The operator $D$ is both a differential and a derivation. 
\end {prop}

\begin{proof}
We show first that $D$ is a differential, that is 
$D^2:E_{\bullet}^n\to E_{\bullet}^{n+2}$ is zero. If $n$ is 
even, 
$D^2=(d +\iota)d=dd+d\iota=2d\iota=0$, the same if $n$ is odd. 
Let's see that $D$ is a derivation, that is 
$D(xy)= D(x)y+ (-1)^{\operatorname{deg}(x)}xD(y).$
There are three cases. 
\begin{enumerate}
 \item Both deg($x$) and deg($y$) even: 
   $$D(xy)=d(xy)=d(x)y + xd(y)=D(x)y + xD(y)$$
 \item  deg($x$) even, deg($y$) odd:
   \begin{eqnarray*}D(xy)&=&(d+\iota)(xy)=d(xy)+\iota xy\\
                        &=& (d(x)y + xd(y))+ \iota xy = d(x)y + (xd(y)+ \iota xy)\\
                        &=&D(x)y + xD(y)
   \end{eqnarray*}
 \item Both deg($x$) and deg($y$) odd:
   \begin{eqnarray*}D(xy)&=&d(xy)=d(x)y + xd(y)=d(x)y + xd(y)+ 2\iota xy\\                                    &=&(d(x)y + \iota xy) + (xd(y)+ \iota xy)\\
                         &=&D(x)y +xD(y).
   \end{eqnarray*}
\end{enumerate} 
For the last case we used the relation $2\iota =0$.
\end{proof}

To prove that the category $\mathcal{W}_A$ has an initial object 
we have to show 
that the solution set condition is satisfied. That means we have to find a 
\emph{set} of objects $\{O_i\}_{i\in I}$, such that for any other object 
$X$ in the category, there is an index $i\in I$ and a map $\phi : O_i \to X$,
not necessarily unique.
Since at each level, a Witt complex $E=E_{\bullet}^*$ is also a complex, there is a map $\lambda :\tilde{\Omega}^* _{W_{\bullet}(A)}\to E_{\bullet}^*$ which in degree zero is the map $\lambda :W_{\bullet}(A)\to E_{\bullet}^0$ prescribed in the definition of a Witt complex. We prove that the image of $\lambda$ is a sub-Witt complex of $E_{\bullet}^*.$ Since the isomorphisms classes of such objects form  a set (they are all quotients of 
$\tilde{\Omega}^* _{W_{\bullet}(A)}$), the solution set condition is satisfied and the 
proposition is proved. 

First of all we have to see that $\iota([1]_n) \in \mbox{Im}(\lambda)$. But this is 
so because 
\begin{eqnarray*}\iota([1]_n) &=& {d(\lambda[-1]_n)\over{\lambda[-1]_n}}
                      = {D(\lambda[-1]_n)\over{\lambda[-1]_n}}\\
            &=& {\lambda d[-1]_n\over{\lambda[-1]_n}} \in \mbox{Im}
(\lambda).
\end{eqnarray*}
 
Since $\iota([1]_n) \in \mbox{Im}(\lambda)$ we see that $\mbox{Im}(\lambda)$ is 
closed under $\iota$. It is also closed under $d$ because it is closed 
under $D$ and $d=D$ or $d=D-\iota$ depending on the degree. It remains 
to see that it is closed under $F$ and $V$.

We start with $F$. The Frobenius operator is multiplicative and each 
element in the image of 
$\lambda$ is of the form $\lambda(a^0da^1\dots  da^n)=
\lambda(a^0)d(\lambda(a^1))\dots d(\lambda a^n)$, so it suffices to show that
$F(\lambda(a))$ and $F(d(\lambda(a))$ are in the image of the canonical map. 
Part of the definition of a Witt complex is that $\lambda F = 
F \lambda$ 
for all $a \in A .$ So $F(\lambda(a)) \in \mbox{Im}(\lambda)$. 
Let us prove that $F(d(\lambda(a))) \in \mbox{Im}(\lambda)$. We use the 
formula:
$$a = [a_0]_n +V([a_1]_{n-1})+\cdots+V^{n-1}([a_{n-1}]_1),$$
which shows that
$$F(d\lambda(a)) = F(d \lambda([a_0]_n)) +F(d \lambda(V([a_1]_{n-1})))
+\cdots+ F(d \lambda(V^{n-1}([a_{n-1}]_1))).$$
Recall from the definition of a Witt complex that
$ Fd\lambda([a]_n) = \lambda([a]_{n-1})d\lambda([a]_{n-1}])$. and that both 
$F$ and $V$ commute with $\lambda$. 
\begin{eqnarray*}
  F(d\lambda(a))&=&  F(d \lambda([a_0]_n))+F(dV \lambda([a_1]_{n-1}))+\cdots+
                       F(d(V^{n-1} \lambda([a_{n-1}]_1))\\
                  &=&  \lambda([a_0]_{n-1})d\lambda([a_0]_{n-1}]) +
                       (d+\iota)(\lambda([a_1]_{n-1})+ \cdots+ \\
                  & &  (d+\iota)V^{n-2} \lambda([a_{n-1}]_1);
\end{eqnarray*}
and this sum clearly is in the image of $\lambda$.

Now we prove that $\mbox{Im}\lambda$ is closed under $V$. 
This follows from the Frobenius reciprocity formula. For example 
\[V(\lambda{a^0}d\lambda(a^1))=V(\lambda(a^0))dV^(\lambda(a^1))-\iota V(\lambda(a^0a^1)).\]
More generally
\[V(\lambda(a^0da^1\dots  da^n))=\sum_{i=1}^n\iota V(a^0a^i)\prod_{1\leq j\neq i\leq n}dV(\lambda(a^j)).\]
This again is in the image of $\lambda.$
\end{proof}

\begin{defn}
The initial object in the category $\mathcal{W}_A$  of Witt complexes over $A$ 
is called the \emph{de Rham Witt complex of $A$} and is denoted 
$W_{\bullet} \Omega^*_A$. 
\end{defn}

\begin{prop}\label{degreezero}For every ring $A$ the following assertions hold:
\begin{itemize} 
\item[(i)] the canonical map $\tilde{\Omega}^* _{W_{\bullet}(A)}\to W_{\bullet}\Omega^*_A$ is surjective,
\item[(ii)]the canonical map $\lambda : W_{\bullet}(A) \to W_{\bullet}\Omega^0_A$ is an isomorphism. 
\end{itemize}
\end{prop}

\begin{proof}
Denote for the moment by $E_{\bullet}^*$ the image of the map $\lambda  :
\tilde{\Omega}^* _{W_{\bullet}(A)}\to W_{\bullet}\Omega^*_A$. It is 
a sub-Witt complex of 
$W_{\bullet}\Omega^*_A$, in particular it is an object of the category 
$\mathcal{W}_A.$ Therefore it admits a unique map from the initial object. We consider the composition: 
$$W_{\bullet}\Omega^*_A\to E_{\bullet}^* \to W_{\bullet}\Omega^*_A.$$
Being an endomorphism of the initial object, it has to be the identity map. 
So the second map is surjective, which amounts to the same thing as the map
$\tilde{\Omega}^* _{W_{\bullet}(A)}\to W_{\bullet}\Omega^*_A$ being surjective. 

In degree zero this means that the map $W_{\bullet}(A)\to W_{\bullet}\Omega^0_A$ is surjective. To prove that it is also injective, we consider the Witt complex $E_{\bullet}^*$ defined by $E_n^0=W_n(A)$ and $E_n^i=0,$ for all $i\geq 1$. As 
$W_{\bullet}\Omega^*_A$ is initial in the category $\mathcal{W}_A$, there 
is a morphism $\mu : W_{\bullet}\Omega^*_A\to E_{\bullet}^*$.
The fact that $\mu$ is a morphism means among other conditions that the diagram
\[\xymatrix{&W_n\Omega_A^0\ar[dd]^{\mu^0}\\
W_n(A)\ar[ur]^{\lambda}\ar[dr]^{\lambda}&\\
&E_n^0}\]
commutes. By the definition of $E_{\bullet}^*$ the corresponding $\lambda$ is the identity morphism, so $\mu^0\circ\lambda=1 $ and it follows that $\lambda$ is injective. 
 \end{proof}

\subsection{The standard filtration}
On every   Witt complex $E_{\bullet}^*$ there is a standard filtration (see also \cite{Ill}):
$$\operatorname{Fil}^s E_n^q=V^sE_{n-s}^q+dV^sE_{n-s}^{q-1}.$$
This filtration can be used to set up inductive arguments when computing de Rham-Witt complexes. The important result that allows this is the following. 
\begin{lemma}\label{exactsequence}
The following sequence is exact: 

$$0\to\operatorname{Fil}^sW_n\Omega_A^q\to W_n\Omega_A^q\xrightarrow{R^{n-s}}W_s\Omega_A^q \to 0$$
\end{lemma}

\begin{proof}
First we show that the composition of the two morphisms is zero. Actually the composition is zero for all Witt complexes. This is so because $R$ commutes with the other operators and any Witt complex is by definition zero in levels zero and below: 
\begin{align*}R^{n-s}(\operatorname{Fil}^sE^q_n)&=R^{n-s}(V^sE^q_{n-s}+dV^sE^{q-1}_{n-s})\\
&=V^sR^{n-s}E^q_{n-s}+dV^sR^{n-s}E^{q-1}_{n-s}\\
&\subset V^sE^i_0+dV^sE^{q-i}_0=0
\end{align*}

Once we know that this composition is zero it follows that $R^{n-s}$ induces a morphism 
\[E_n^q/{\operatorname{Fil}^s(E_n^q)}\longrightarrow{R^{n-s}} E_x^q .\]
To end the proof of the lemma we need to show that this morphism is an isomorphism for $E_{\bullet}^*=W_{\bullet}\Omega^*$. Fix a value of $n-s$ and define 
\[W_s'\Omega_A^i=W_n\Omega_A^q/{\operatorname{Fil}^sW_n\Omega^q_A}.\]
We prove that this is a Witt complex over $A$. We only need to check that the operators are well defined, then the relations are automatically satisfied.
To show that $R$ and $F$ induce operators $R,F:W_s'\Omega_A^q\to W_{s-1}'\Omega_A^q$ we need to show $R(\operatorname{Fil}^sW_n\Omega^q_A)\subset\operatorname{Fil}^{s-1}W_{n-1}\Omega^q_A$ and $F(\operatorname{Fil}^sW_n\Omega^q_A)\subset\operatorname{Fil}^{s-1}W_{n-1}\Omega^q_A$. The first relation follows from $VR=RV$ and $dR=Rd$ and the second from $FV=2$ and $FdV=d+\iota$. Similarly $V$ induces an operator on $W_{\bullet}'\Omega_A^*$ if $V(\operatorname{Fil}^sW_n\Omega^q_A)\subset\operatorname{Fil}^{s+1}W_{n+1}\Omega^q_A$ and this follows from $Vd=2dV$. $\iota$ and $d$ induce operators if $\iota(\operatorname{Fil}^sW_n\Omega^q_A)\subset\operatorname{Fil}^{s}W_{n}\Omega_A^{q+1}$ and $d(\operatorname{Fil}^sW_n\Omega^q_A)\subset\operatorname{Fil}^{s}W_{n}\Omega^{q+1}_A$. The first follows from $\iota V=V\iota$ and $\iota d=d\iota$ and the second from $dd=d\iota$.

We show now that $W_{\bullet}'\Omega_R^*$ is an initial object in the category of 2-typical Witt complexes over $A$ and hence the morphism induced by $R^{n-s}$, $W_{\bullet}'\Omega_R^*\to W_{\bullet}\Omega_R^*$ is an isomorphism and hence the sequence is exact. 

Consider $E_{\bullet}^*$ a 2-typical Witt complex over $A$. We construct a morphism $W_{\bullet}'\Omega_A^*\to E_{\bullet}^*$ and show that it is unique. Since the standard filtration is natural we have maps: 
\[W_s'\Omega^i_A=W_n\Omega^i_A/{\operatorname{Fil}^sW_n\Omega^q_A}\to E^i_n/{\operatorname{Fil}^sE^i_n}\xrightarrow{R^{n-s}}E^i_s. \]

To show that this homomorphism of Witt complexes is unique we first show that the map $\tilde{\Omega}^i(W_s(A))\to W_s'\Omega^i_A$ is surjective. This is immediate from the diagram
\[\xymatrix{\tilde{\Omega}^i(W_n(A))\ar@{>>}[d]\ar@{>>}[r]&W_n\Omega_A^i\ar@{>>}[d]\\\tilde{\Omega}^i(W_s(A))\ar[r]&W_s'\Omega_A^i}\]
Now in the diagram
\[\xymatrix{\tilde{\Omega}^*(W_{\bullet}(A))\ar[r]\ar@{>>}[d]&E_{\bullet}^*\\W_{\bullet}'\Omega^*_A\ar[ur]&}\]
consider in the category of pro-DGA the top map is unique, therefore the oblique map is unique. 
\end{proof}

\subsection{An additivity result} As far as we know, the relations in the definition of a 2-typical Witt complex are independent. However, one relation can be partially deduced from the others. We make this precise in the following Lemma. 

\begin{lemma}\label{additivity}
Let $A$ be an arbitrary ring, and $E_\bullet^*$ a pro-graded ring. Assume that $E_\bullet^*$ is endowed with all the operators in the definition of a 2-typical Witt complex, and all the relations are sastisfied with the exception of the last relation. Assume that this relation holds for two given elements $f, g\in A,$ that is
\[Fd\lambda([f]_n) = \lambda([f]_{n-1})d\lambda([f]_{n-1}]),\]
\[Fd\lambda([g]_n) = \lambda([g]_{n-1})d\lambda([g]_{n-1}]).\]
Then it also holds for their sum, $f+g\in A:$
\[Fd\lambda([f+g]_n) = \lambda([f+g]_{n-1})d\lambda([f+g]_{n-1}]).\]
\end{lemma}
\begin{proof} 
The proof is inspired by the proof of Proposition 1.3 in \cite{LZ}.

Since there is no danger of confusion, we omit $\lambda.$

We prove the statement by induction on the level $n.$ If $n=1$ the relation holds trivially. Assume we have proved the Lemma for $n-1.$ We know that
\[Fd[f]_n = [f]_{n-1}d[f]_{n-1}],\]
\[Fd[g]_n = [g]_{n-1}d[g]_{n-1}],\]
and if we apply $R$ to these relations we obtain
\[Fd[f]_{n-1} = [f]_{n-2}d[f]_{n-2}],\]
\[Fd[g]_{n-1} = [g]_{n-2}d[g]_{n-2}].\]
By the induction hypothesis, we have
\[Fd[f+g]_{n-1}=[f+g]_{n-2}d[f+g]_{n-2}.\]
We define $\tau\in W_{n-1}(A)$ by the formula: $$[f+g]_n=[f]_n+[g]_n+V\tau.$$ We apply $R$ to both sides of this identity: $$[f+g]_{n-1}=[f]_{n-1}+[g]_{n-1}+VR\tau.$$
We square both sides: 
\[[f+g]_{n-1}^2=[f]_{n-1}^2+[g]_{n-1}^2+(VR\tau)^2+2([f]_{n-1}+[g]_{n-1})VR\tau+2[fg]_{n-1},\]
and, since $(VR\tau)^2=V(R\tau FVR\tau)=2VR(\tau^2)$, we get:
\[F([f+g]_n+[f]_n-[g]_n)=2(VR\tau +[fg]_{n-1}+([f]_{n-1}+[g]_{n-1})VR\tau).\]
The left hand side of this identity is $F(V\tau)=2\tau$, thus we obtain: 
\[2\tau=2(VR\tau +[fg]_{n-1}+([f]_{n-1}+[g]_{n-1})VR\tau).\]
This identity is true in the ring $W_{n-1}(R)$ for any ring $R$, in particular for the ring of polynomials in two variables $R=\Z[f,g]$. For this ring the multiplication by 2 in $W_{n-1}(R)$ is injective, therefore the following identity is true for this particular ring:   
\[\tau=VR\tau +[fg]_{n-1}+([f]_{n-1}+[g]_{n-1})VR\tau.\]
Taking Witt vectors of length $n-1$ is functorial, and it follows that the identity is true for arbitrary rings $R$ and elements $f$ and $g$.

Once we proved this formula we return to the identity we want to prove: 
$$Fd[f+g]_n=[f+g]_{n-1}d[f+g]_{n-1},$$
or equivalently:
$$Fd([f]_n+[g]_n+V\tau)=([f]_{n-1}+[g]_{n-1}+VR\tau)(d[f]_{n-1}+d[g]_{n-1}+dVR\tau).$$
We expand the right hand side:
\begin{align*}
Fd[f]_n+Fd[g]_n+FdV\tau&=[f]_{n-1}d[f]_{n-1}+[g]_{n-1}d[g]_{n-1}\\
&+d([fg]_{n-1}+([f]_{n-1}+[g]_{n-1})VR\tau)+VR\tau dVR\tau.
\end{align*}
Using the hypothesis that $Fd[f]_n=[f]_{n-1}d[f]_{n-1}$ and $Fd[g]_n=[g]_{n-1}d[g]_{n-1}$, and the formula for $\tau$, the previous identity becomes equivalent to:
$$FdV\tau=d(\tau-VR\tau^2)+VR\tau dVR\tau$$
or:
$$d\tau+\iota\tau=d\tau-dVR\tau^2+VR\tau dVR\tau.$$
We reduced the problem to proving the following formula: 
$$\iota\tau=-dVR\tau^2+VR\tau dVR\tau.$$
We prove this separately in the next Lemma.
\end{proof}
\begin{lemma}For every ring $A$ and every element $\tau\in W_{k}(A)$, the following identity holds: 
\[dVR\tau^2=VR\tau dVR\tau+\iota\tau\]
\end{lemma}
\textbf{Remark:} This lemma says that $\iota$ measures the failure of $d$ to be a PD-derivation. To explain this we need to recall what a PD-structure on a ring is and what a PD-derivation is. 

If $A$ is a commutative and unital ring and $I$ is and ideal, a PD-structure on $(A,I)$ is a family of maps $\gamma_n:I\to A$, $n\in \N$ which morally behave like dividing the $n$'th power by $n!$. More precisely they are required to satisfy five conditions: 
\begin{itemize}
\item[(i)] $\gamma_0(x)=1,$ $\gamma_1(x)=x,$ $\gamma_n(x)\in I,$ $\forall x\in I,$
\item[(ii)]$\gamma_n(x+y)=\sum_{i=0}^n\gamma_i(x)\gamma_{n-i}(y), \forall x,y\in I,$
\item[(iii)]$\gamma_n(ax)=a^n\gamma_n(x), \forall a\in A, x\in I,$
\item[(iv)]$\gamma_p(x)\gamma_q(x)=\binom{p+q}{p}\gamma_{p+q}(x),\forall p,q\in\N, x\in I,$
\item[(v)]$\gamma_p(\gamma_q(x))=\frac{(pq)!}{p!(q!)^p}\gamma_{pq}(x),\forall p,q\in\N, x\in I.$
\end{itemize}

If $A$ is a $\Z_{(p)}$-algebra, there is a canonical PD-structure on $(W_n(A),VW_{n-1}(A)),$ namely:
\[\gamma_m:VW_{n-1}(A)\to W_n(A)\]
\[\gamma_m(Vx)=\begin{cases}
                1,& \text{if } m=0,\\
                \frac{p^{m-1}}{m!}V(x^m),&\text{if } m\geq 1
               \end{cases}\]

If $(A,I)$ is a ring with a PD-structure and $d:A\to M$ is a derivation of $A$ into an $A$-module $M$, then $d$ is called a PD-derivation if  $d(\gamma_n(x))=\gamma_{n-1}(x)dx,$ for all $x\in I.$
If we consider $W(A),$ the inverse limit of the pro-ring $W_\bullet(A)$, elements in it are sequences of elements $x_n\in W_n(A)$, and$R: W(A)\to W(A)$ becomes the identity morphism. So the identity in the lemma reads: if $x\in VW(A)\subseteq W(A),$ $x=V\tau,$ then:
\[d(\gamma_2(x))=\gamma_1(x)dx+\iota\tau.\]
If we did not have $\iota\tau$ in this relation, then  Lemma 1.2 of Langer, Zink \cite{LZ}  would show that $d$ is a PD-derivation. 
We will now prove the lemma. 
\begin{proof}
The proof is in three steps. First we show that the identity holds for all elements $\tau=[\phi]_{k}$, for $\phi\in A[X]$. Then we show that once it holds for $\tau$ it also holds for $V\tau.$ Finally we show that if the identity holds for  two elements $\tau_1$, $\tau_2$ it also holds for their sum $\tau_1+\tau_2.$

We begin with $\tau=[\phi]_{k}$ for some $\phi\in A[X]$. We manipulate the left hand side of the identity that we want to prove: 
\begin{align*}
dV[\phi]^2_{k-1}&=dV(F([\phi]_k))=d(V([1]_{k-1})[\phi]_k)\\
&=dV([1]_{k-1})[\phi]_k+V([1]_{k-1})d[\phi]_{k}\\
&=dV([1]_{k-1})[\phi]_k+V(Fd[\phi]_{k}),
\intertext{and using the induction hypothesis that $Fd[\phi]_{k}=[\phi]_{k-1}d[\phi]_{k-1}$ for $k\leq n-1$ :}
dV[\phi]^2_{k-1}&=(dV[1]_{k-1})[\phi]_k+V([\phi]_{k-1}d[\phi]_{k-1}).
\intertext{The right hand side is} 
V[\phi]_{k-1}dV[\phi]_{k-1}+\iota[\phi]_{k}&=V([\phi]_{k-1}FdV[\phi]_{k-1})+\iota[\phi]_{k}\\
&=V([\phi]_{k-1}d[\phi]_{k-1})+V(\iota[\phi]_{n-1}^2)+\iota[\phi]_k.
\intertext{Therefore the equality of the two members is equivalent to the equality}
dV[1]_{k-1}[\phi]_k&=V(\iota[1]_{k-1})[\phi]_n+\iota[\phi]_k,
\intertext{which is certainly true if}
dV[1]_{k-1}&=V(\iota[1]_{k-1})+\iota[1]_k.
\end{align*}

This last identity follows from Lemma \ref{iota} which gives the formula for $\iota[1]_k$.

Now we assume we know the relation for $\tau$ and we want to prove it for $V(\tau)$. We know:
\begin{align*}
dVR\tau^2&=VR\tau dVR\tau+\iota\tau,
\intertext{we apply $V$ to this:}
V(dVR\tau^2)&=V(VR\tau dVR\tau)+V(\iota\tau);
\intertext{The left hand side of this equation is:}
V(dVR\tau^2)&=2dV^2R\tau^2=dV^2(R\tau FVR\tau)\\
&=dV(V(R\tau)V(R\tau))=dVR(V(\tau)^2)
\intertext{We want to prove:}
dVR(V(\tau)^2)&=VRV\tau dVRV\tau+\iota V\tau.
\intertext{We notice that the left hand side of this equation is equal to the left hand side of the previous equation, so it is sufficient for us to prove:}
V(VR\tau dVR\tau)+V(\iota\tau)&=VRV\tau dVRV\tau+\iota V\tau,
\intertext{or:}V(VR\tau dVR\tau)&=VRV\tau dVRV\tau\\
V(VR\tau dVR\tau)&=V((RV\tau)FdVRV\tau)\\
V(VR\tau dVR\tau)&=V((RV\tau)dRV\tau)+V((RV\tau)\iota(RV\tau))\\
0&=\iota(V((VR\tau)^2))\\
0&=\iota(V(R\tau FVR\tau))
\intertext{which is true since $FV=2$ and $2\iota=0.$}
\end{align*}

Finally we want to prove that if the relation holds for $\tau_1$ and $\tau_2$, it also holds for $\tau_1+\tau_2.$
We know:
\begin{align*}
RdV\tau_1^2&=RV\tau_1dRV\tau_1+\iota\tau_1,\\
RdV\tau_2^2&=RV\tau_2dRV\tau_2+\iota\tau_2
\intertext{and we want:}
RdV(\tau_1+\tau_2)^2&=RV(\tau_1+\tau_2)dRV(\tau_1+\tau_2)+\iota(\tau_1+\tau_2),
\intertext{or equivalently:}
RdV\tau_1^2+RdV\tau_2^2+2RdV(\tau_1\tau_2)&=RV\tau_1dRV\tau_1+RV\tau_2dRV\tau_2\\
&\phantom{=}+d(RV\tau_1RV\tau_2)+\iota\tau_1+\iota\tau_2.
\intertext{After cancelling six terms the equation reduces to:}
2RdV\tau_1\tau_2&=RdV\tau_1V\tau_2,\\
2RdV\tau_1\tau_2&=RdV(\tau_1FV\tau_2),
\end{align*}
which is true since $FV=2$.
\end{proof}

\section{Computations}
\subsection{The 2-typical de Rham-Witt vectors of the integers}

Before we state the structure theorem for 
$W_{\bullet}\Omega^*_\field{Z}$ we introduce a bit of notation. 
We denote by $V^i(1)\in W_n\Omega^0_A$ the element $V^i([1]_{n-i})$.

\begin{thm}\label{structureforZ}
The structure of $W_{\bullet}\Omega^*_\field{Z}$ is as follows 

\begin{itemize}
\item[(i)]As abelian groups
\begin{eqnarray}
W_n\Omega^0_\field{Z}&=& \bigoplus_{i=0}^{n-1}\field{Z}\cdot V^i(1),\\
W_n\Omega^1_\field{Z}&=& \bigoplus_{i=1}^{n-1}\field{Z}/{2^i\field{Z}}\cdot dV^i(1) ,\\
W_n\Omega^i_\field{Z}&=&0,\quad \mbox{for } i\geq 2.
\end{eqnarray}
\item[(ii)]The product is given by 
\begin{eqnarray}
V^i(1)\cdot V^j(1)& =&2^iV^j(1), \quad\mbox{if}\quad i \leq j,\\
V^i(1)\cdot dV^j(1)& =& 
\begin{cases} 2^idV^j(1)+\sum_{s=j+1}^{n-1}2^{s-1}dV^s(1),&\text{if $1\leq i< j$}\\
                    \sum_{s=i+1}^{n-1}2^{s-1}dV^{s}(1) ,&\text{if $1\leq j\leq i$}.\\
  \end{cases}\\
\end{eqnarray}
\item[(iii)]The operator $V$ acts as follows
\begin{eqnarray}
V(V^i(1))&=&V^{i+1}(1),\\
V(dV^i(1))&=&2dV^{i+1}(1).
\end{eqnarray}
\item[(iv)]The operator $F$ acts as follows
\begin{eqnarray}
F(V^i(1))&=&2V^{i-1}(1),\\
F(dV^i(1))&=&dV^{i-1}(1)+\sum_{s=i}^{n-2}2^{s-1}dV^{s}(1).
\end{eqnarray}
\item[(v)]The operator $d$ acts by $d(V^i(1))=dV^i(1)$ when $i\geq 1$ and $d(1)=0$, and the action of the operator $\iota$ is given by
\begin{eqnarray}
\iota(V^i(1))&=&\sum_{s=i+1}^{n-1}2^{s-1}dV^{s}(1).
\end{eqnarray}
\item[(vi)]The operator $R: W_n\Omega^*_\field{Z}\to W_{n-1}\Omega^*_\field{Z}$
acts as follows
\begin{eqnarray}
RV^i(1)&=&\begin{cases}V^i(1)&\text{if $i\leq n-2$}\\
        0&\text{if $i=n-1$}\\
        \end{cases}\\
RdV^i(1)&=&\begin{cases}dV^i(1)&\text{if $i\leq n-2$}\\
        0&\text{if $i=n-1$.}\\
        \end{cases}
\end{eqnarray}

\end{itemize}
\end{thm}

\begin{proof}

We begin with the fifth assertion and  then we prove the others in the stated order. 

(v) We already know the formula (see \ref{iota}):
$$\iota[1]_n=\sum_{s=1}^{n-1}2^{s-1}dV^{s}(1),$$
which is the particular case for the relation we want to prove when $i=0$. 
For other $i$ we have:
\begin{eqnarray*}
\iota V^i([1]_{n-i})&=&V^i(\iota[1]_{n-i})\\
                    &=&V^i(\sum_{s=1}^{n-i-1}2^{s-1}dV^{s}([1]_{n-i-s-1})\\
                    &=&\sum_{s=1}^{n-i-1}2^{s+i-1}dV^{s+i}([1]_{n-i-s-1})\\  
                    &=&\sum_{s=i+1}^{n-2}2^{s-1}dV^{s}([1]_{n-s-1}).
\end{eqnarray*}
(i),(ii)The isomorphism described in the first relation follows from the previous theorem. The second relation follows from the fact that the map $\lambda: \tilde{\Omega}^* _{W_{\bullet}(A)}\to W_{\bullet}\Omega^*_A$ is surjective and from the product relations that we now prove. 
The first product relation is the product rule described in \ref{wittbasis} in the case $p=2$. 
The second product relation: 
\begin{itemize}
\item If $1\leq i \leq j$:
\begin{eqnarray*}
V^i(1)\cdot dV^j(1)&=&V^i(F^id V^j(1))=V^i((d+\iota)V^{j-i}(1))\\
                   &=&V^idV^{j-i}(1)+V^j(\iota (1))\\
                   &=&2^idV^j(1)+V^j(\sum_{s=1}^{n-j-1}2^{s-1}dV^{s}(1))\\
                   &=&2^idV^j(1)+\sum_{s=1}^{n-j-1}2^{s+j-1}dV^{j+s}(1)\\
                   &=&2^idV^j(1)+\sum_{s=j+1}^{n-1}2^{s-1}dV^{s}(1).
\end{eqnarray*}
\item If $1\leq j\leq i$:

\begin{eqnarray*}
V^i(1)\cdot dV^j(1)&=&d(V^i(1)V^j(1))-dV^i(1)V^j(1)\\
                   &=&d(2^jV^i(1))-V^j(1)dV^i(1)\\
                   &=&2^jdV^i(1)-2^jdV^i(1)-\sum_{s=i+1}^{n-1}2^{s-1}dV^{s}(1)\\
                   &=&\sum_{s=i+1}^{n-1}2^{s-1}dV^{s}(1).
\end{eqnarray*}
\end{itemize} 
 
We note here that one can give a unified product relation for $V^i(1)\cdot dV^j(1)$, namely: 
$$V^i(1)\cdot dV^j(1)=2^idV^j(1)+\sum_{s=\max(i,j)+1}^{n-1}2^{s-1}dV^{s}(1)$$

(iii) The first relation is trivial and the second follows from $Vd=2dV$.

(iv) The first relation follows from $FV=2$. The second relation:
\begin{eqnarray*}
F(dV^i(1))&=&(d+\iota)V^{i-1}(1)\\
          &=&dV^{i-1}(1)+\iota(V^{i-1}(1)\\
          &=&dV^{i-1}(1)+\sum_{s=i-1}^{n-2}2^{s}dV^{s+1}(1).
\end{eqnarray*}

Once we have these relations and the fact that $2^idV^i(1)=V^id(1)=0$, it follows that $\lambda$ factors through a surjective
map $\bigoplus_{i=1}^{n-1}\field{Z}/{2^i\field{Z}}\cdot dV^i(1)\to W_{\bullet}\Omega^1_\field{Z}$. 

We prove now that $W_{\bullet}\Omega^i_\field{Z} = 0$, for $i\geq 2$. We will prove this by induction on the level using the 
standard filtration.
The first step of the induction, that 
$W_1\Omega^q_\field{Z}=0,$ forall $q\geq 2$ follows from the surjectivity of 
the map $\lambda: \Omega^q _{W_1(\field{Z})}=\Omega^q _\field{Z}\to W_1\Omega^q_\field{Z}$, and 
fact that the domain of the map is zero whenever $q\geq 1$. Assuming that   
$W_n\Omega^q_\field{Z}=0$ for all $q\geq 2$ we prove that $W_{n+1}\Omega^q_\field{Z}=0,$ for all $q\geq 2$. This is so because in the short exact sequence
$$0\to\operatorname{Fil}^{n-1}W_n\Omega_\field{Z}^q\to W_n\Omega_\field{Z}^q\xrightarrow{R}W_{n-1}\Omega_\field{Z}^q \to 0$$
the right term is zero by induction and the left term is zero because
$\operatorname{Fil}^{n-1}W_n\Omega_\field{Z}^q=V^{n-1}W_1\Omega_\field{Z}^2+dV^{n-1}W_1\Omega_\field{Z}^1$ and both $W_1\Omega_\field{Z}^1$ and $W_1\Omega_\field{Z}^2$
 are zero as seen above. 

To finish the description of the groups that form the de Rham-Witt complex 
$W_{\bullet}\Omega^*_\field{Z}$ we consider the pro-graded ring 
$G_{\bullet}^*$ defined by the groups on the right hand side of
the relations $(1)-(3)$, that is: 
\begin{eqnarray*}
G_n^0&=& \bigoplus_{i=0}^{n-1}\field{Z}\cdot V^i(1),\\
G_N^1&=& \bigoplus_{i=1}^{n-1}\field{Z}/{2^i\field{Z}}\cdot dV^i(1) ,\\
G_n^i&=&0,\quad \mbox{for } i\geq 2.
\end{eqnarray*}
The product is defined by the relations in (ii) 
, the operators $F, V, d, \iota, R$ are given by the relations in 
(iii)- (vi). We check that with these definitions $G_{\bullet}^*$ 
is indeed a Witt complex. 

The only non-trivial relation to verify  is that $ Fd\lambda([a]_n) = \lambda([a]_{n-1})d\lambda([a]_{n-1}]),$ for all integers $a$. Using the additivity result \ref{additivity}, we see that we need to check this relation only for the integers $a=1$ and $a=-1.$ It is trivially satisfied in the first case, and easy to see in the second, once we recall from \ref{relation-1} that $[-1]_n= -[1]_n +V[1]_{n-1}:$ 
\begin {align*}
Fd[-1]_n&=Fd(-[1]_n+V[1]_{n-1})=FdV[1]_{n-1}\\
&=(d+\iota)[1]_{n-1}=\iota[1]_{n-1}\\
&=[-1]_{n-1}d[-1]_{n-1}.
\end{align*}

To prove now that $W_{\bullet}\Omega^*_\field{Z}\cong G_{\bullet}^*$ we define a morphism of Witt complexes 
\begin{align*}
G_{\bullet}^*&\longrightarrow W_{\bullet}\Omega^*_\field{Z}\\
 V^i(1)                 &\longmapsto V^i(1) \\
dV^i(1)                 &\longmapsto dV^i(1)
\end{align*}
The composition $W_{\bullet}\Omega^*_\field{Z}\to G_{\bullet}^*\to
W_{\bullet}\Omega^*_\field{Z}$ is an endomorphism of the initial object in the category $\mathcal{W}_A$ and as so it is the identity. The composition $G_{\bullet}^*\to W_{\bullet}\Omega^*_\field{Z}\to G_{\bullet}^*$
is an endomorphism of $G_{\bullet}^*$; it is not hard to see that  
the only endomorphism of $G_{\bullet}^*$ is the identity: being a 
morphism of pro-rings it maps $[1]_n$ to itself, and since it commutes with $V$ and $d$ it will also map $V^i(1)$ and $dV^i(1)$ to themselves.  
\end{proof}

\textbf{Remark:} The same proof works to give us the structure of $W_{\bullet}\Omega^*_{\Z_{(2)}}:$
\begin{itemize}
\item[(i)]As abelian groups: 
\begin{eqnarray}
W_n\Omega^0_{\Z_{(2)}}&=& \bigoplus_{i=0}^{n-1}\Z_{(2)}\cdot V^i(1),\\
W_n\Omega^1_{\Z_{(2)}}&=& \bigoplus_{i=1}^{n-1}\Z/{2^i\Z}\cdot dV^i(1) ,\\
W_n\Omega^i_{\Z{(2)}}&=&0,\quad \mbox{for } i\geq 2.
\end{eqnarray}
\item[(ii)] The product formulas and the actions of the various operators are the same as in Theorem \ref{structureforZ}.
\end{itemize}

Indeed, the only thing we have to check is that $\Omega^i_{\Z_{(2)}}=0$ for all $i\geq 2.$ To see this we need to prove that $d(\frac{1}{m})=0$ for all $m\in \Z$ odd. This follows from $0=d(1)=d(m\frac{1}{m})=md(\frac{1}{m}),$ 
since $m\in \Z$ is a unit.

\subsection{The 2-typical de Rham-Witt complex for polynomial extensions}
In this subsection we describe the relationship between the 2-typical de Rham-Witt complex of the ring of polynomials in one variable over a $\Z_{(2)}$-algebra $A$ and the 2-typical de Rham-Witt complex of the ring $A$. In order to do that we will identify the left adjoint of the forgetful functor $\mathcal{W}_{A[X]}\to\mathcal{W}_A$. We call this functor $P:\mathcal{W}_A\to\mathcal{W}_{A[X]}$, and since it commutes with colimits, it will carry $W_{\bullet}\Omega^*_A$ into $W_{\bullet}\Omega^*_{A[X]}$. 

In order to define the functor $P$ we first analyze the Witt complex $W_{\bullet}\Omega^*_{A[X]}$. Inside it we find the image of the map $W_{\bullet}\Omega^*_{A}\to W_{\bullet}\Omega^*_{A[X]}$ induced by the inclusion $A\to A[X]$. Besides this image we can certainly identify the elements $[X]^i_n$. If we play with the multiplication and with the operators $R$, $F$, $V$, $d$, and $\iota$ we will find new elements, but because of the relations that hold in every Witt complex, we will see that all these elements can be classified in four types. The first obvious type is the elements of the form $a[X]^i_n$, where $a\in \operatorname{Im}(W_{\bullet}\Omega^*_{A})$, $i\in \N$, and $n \geq 1$. When there is no danger of confusion, we  omit the subscript $n$, and also write $a\in W_{\bullet}\Omega^*_{A}$. This type is closed under multiplication and also under the action of $R$, $F$, and $\iota$. If we apply $d$ and $V$ we will get two new types: elements of the form $b[X]^{k-1}d[X]$ and elements of the form $V^r(c[X]^l)$, where $b, c \in W_{\bullet}\Omega^*_{A}$, and $k, r, l >0$. Special attention has to be paid to the latter type, as some elements of that form were already listed as elements of the first type. An example is $V(c[X]^2) = V(cF([X]))= V(c)[X]$. The restriction that we have to impose is $l$ be odd. Finally, if we apply  $V$ and then $d$ we obtain a new type, of elements of the form $dV^s(e[X]^m)$, where $s>0$, $e\in W_{\bullet}\Omega^*_{A}$, and $m$ is odd. If we multiply elements of any of these two types together we will get a sum of elemtents of these types. The key observation is the following:

\begin{lemma}
In any Witt complex over $A[X]$ the following relation holds:
$$d[X]d[X]=\iota([1])[X]d[X].$$
\end{lemma}
\begin{proof}
\begin{align*}
d[X]d[X]&=d([X]d[X])-[X]dd[X]=d(Fd[X])-[X]d\iota[X]\\
&=2Fdd[X]+\iota([1])[X]d[X]=2Fd\iota[X]+\iota([1])[X]d[X]\\
&=\iota([1])[X]d[X].
\end{align*}
\end{proof}
Using this observation we can see for example that the product of two elements of the second type is again an element of second type: 
\[b[X]^{k-1}d[X]b'[X]^{k'-1}d[X]=\iota(bb')[X]^{k+k'-1}d[X].\] 
The other products and the action of the different operators on the elements can also be derived. The formulas that we obtain will be exactly the formulas that we plug in the definition of the functor $P$. 

Before we define the functor $P$ we need to recall a result of Hesselholt, Madsen that describes the ring of $p$-typical Witt vectors over the ring $A[X]$: every element $f\in W_n(A[X])$ can be written uniquely as a sum:
\[f=\sum_{j\in \N}a_{0,j}[X]^j_n+\sum_{s=1}^{n-1}\sum_{(j,p)=1}V^s(a_{s,j}[X]_{n-s}^j),\] with $a_{s,j}\in W_{n-s(A)}$ and all but finitely many $a_{s,j}$ zero (see Lemma 4.1.1 in \cite{He1}). In the case $p=2$ this results read: every element $f\in W_n(A[X])$ can be written uniquely as a finite sum of elements of two types, that we will call type 1 and type 3, for reasons that will soon become clear: 
 
 \begin{itemize}
\item Type 1: elements of the form $a[X]^j$, where $a\in W_n\Omega^q$,
\item Type 3: elements of the form $V^r(c[X]^l)$, where $r>0$, $c\in W_{n-r}\Omega^q$, and $l$ is odd.
 \end{itemize}
Now we are ready to define the functor $P:\mathcal{W}_A\to\mathcal{W}_{A[X]}$. On objects it is defined as follows: for a Witt complex $E^*_{\bullet}\in \mathcal{W}_A$, $P(E)_n^q$ consists of formal sums of four types of elements: 
 \begin{itemize}
\item Type 1: elements of the form $a[X]^j$, where $a\in E_n^q$,
\item Type 2: elements of the form $b[X]^{k-1}d[X]$, where $b\in E_n^{q-1}$,
\item Type 3: elements of the form $V^r(c[X]^l)$, where $r>0$, $c\in E_{n-r}^q$, and $l$ is odd,
\item Type 4: elements of the form $dV^s(e[X]^m)$, where $s>0$, $e\in E_{n-s}^{q-1}$, and $m$ is odd.
 \end{itemize}
The product is graded commutative, and is given by the following ten formulas: 
\begin{itemize}
\item[P1.1:] $a[X]^ja'[X]^{j'}=aa'[X]^{j+j'}$,
\item[P1.2:] $a[X]^jb[X]^{k-1}d[X]=ab[X]^{j+k-1}d[X]$,
\item[P1.3:] $a[X]^jV^r(c[X]^l)=V^r(F^r(a)c[X]^{2^rj+l})$,
\item[P1.4:] 
\begin{align*}
a[X]^jdV^s(e[X]^m)=&(-1)^{|a|}\frac{m}{2^sj+m}dV^s(F^s(a)e[X]^{2^sj+m})\\
&-(-1)^{|a|}V^s((F^s(da)e-\frac{j}{2^sj+m}d(F^s(a)e))[X]^{2^sj+m}),
\end{align*}
\item[P2.2:]$b[X]^{k-1}d[X]b'[X]^{k'-1}d[X]=\iota(bb')[X]^{k+k'-1}d[X]$,
\item[P2.3:]
\begin{align*}
b[X]^{k-1}d[X]V^r(c[X]^l)=&-(-1)^{|b|}\frac{1}{2^rk+l}V^r(d(F^r(b)c)[X]^{2^rk+l})\\
&+(-1)^{|b|}\frac{2^r}{2^rk+l}dV^r(F^r(b)c[X]^{2^rk+l}),
\end{align*}
\item[P2.4:]
\begin{align*}
b[X]^{k-1}d[X]dV^s(e[X]^m)=&-(-1)^{|b|}\frac{1}{2^sk+m}V^s(F^s(db+k\iota(b))de[X]^{2^sk+m})\\
&+(-1)^{|b|}\frac{1}{2^sk+m}dV^s(F^s(b)de[X]^{2^sk+m}),
\end{align*}
\item[P3.3:]
\begin{align*}
\hspace{1cm}\lefteqn{V^r(c[X]^l)V^{r'}(c'[X]^{l'})=}\\
&=\begin{cases}
          2^{r'}V^r(cF^{r-r'}(c')[X]^{2^{r-r'}l'+l})   ,&\text{if $r>r'$,} \\
          2^rV^{r-v}(V^v(cc')[X]^{2^{-v}(l+l')})   ,& \text{if $r=r'$ and $v=v_2(l+l'), v\leq r$,}\\
          2^rV^r(cc')[X]^{2^{-r(l+l')}}, &\text{if $r=r'$ and $v>r$,}
   \end{cases}
\end{align*}
\item[P3.4:]
\begin{align*}
\hspace{1cm}\lefteqn{V^r(c[X]^l)dV^s(e[X]^m)=}\\
&\hspace{-1.5cm}=\begin{cases}
\lefteqn{(-1)^{|c|}\frac{2^rm}{2^{s-r}l+m}dV^s(F^{s-r}(c)e[X]^{2^{s-r}l+m}),} & \hspace{1cm}\text{if $r<s$,}\\
\lefteqn{V^{r-v}(V^v(c(d+\iota)(e))[X]^{2^{-v}(l+m)})}&\\
\hspace{.5cm}+(-1)^{|c|}\frac{2^rm}{l+m}dV^r(ce[X]^{2^{-v}(l+m})&\\
\hspace{.5cm}-(-1)^{|c|}\frac{2^vm}{l+m}V^{r-v}(dV^v(ce)[X]^{2^{-v}(l+m})),& \text{if $r=s$, $v=v_2(l+m)<r$,}\\
\lefteqn{V^r((c(d+\iota)(e))[X]^{2^{-r}(l+m)})}&\\
\hspace{.5cm}+(-1)^{|e|}mV^r(ce)[X]^{2^{-r}(l+m)-1}d[X],&\text{if $r=s$, $v=v_2(l+m)\geq r$,}\\
\lefteqn{V^r(cF^{r-s}((d+\iota)(e))[X]^{2^{r-s}m+l})}&\\
\hspace{.5cm}\lefteqn{+(-1)^{|c|}\frac{2^rm}{2^{r-s}m+l}dV^r(cF^{r-s}(e)[X]^{2^{r-s}m+l})}&\\
\hspace{.5cm}\lefteqn{-(-1)^{|c|}\frac{m}{2^{r-s}m+l}V^r(d(cF^{r-s}(e))[X]^{2^{r-s}m+l}),}&\hspace{2cm}\text{if $r>s$,}
\end{cases}
\end{align*}
\item[P4.4:]
\begin{align*}
\hspace{1cm}\lefteqn{dV^s(e[X]^m)dV^{s'}(e'[X]^{m'})=}\\
&\hspace{-1.5cm}=\begin{cases}
\lefteqn{-(-1)^{|e|}dV^{s'}((F^{s'-s}((d+\iota)(e)e')+\frac{m}{2^{s'-s}m+m'}d(F^{s'-s}(e)e'))[X]^{2^{s'-s}m+m'})}&\\
+V^{s'}((F^{s'-s}(de)\iota(e')+F^{s'-s}(e)d\iota(e'))[X]^{2^{s'-s}m+m'}),&\text{if $s<s'$,} \\
V^{s-v}((V^v(ed\iota(e'))+dV^v(e\iota(e')))[X]^{2^{-v}(m+m')})&\\
+dV^{s-v}((V^v(e(d+\iota)(e'))+dV^v(ee'))[X]^{2^{-v}(m+m')}),& \text{if $s=s'$, $v=v_2(l+m)<s$,}\\
\lefteqn{(V^s(e\iota(e'))+(-1)^{e'}m'dV^s(ee')[X]^{2^{-s}(m+m')-1}d[X])}&\\
+V^s(ed\iota(e')[X]^{2^{-s}(m+m')})&\\
+dV^s(e(d+\iota)(e')[X]^{2^{-s}(m+m')}),& \text{if $s=s'$, $v=v_2(l+m)\geq s$.}\\
\end{cases}
\end{align*}
\end{itemize}
The definition of $\lambda$, $R$, and $\iota$ are obvious, the action of $V$ is given by the following four formulas: 
\begin{itemize}
\item[V1:] \[V(a[X]^j)=
  \begin{cases}V(a[X]^j),& \text{if $j$ odd,}\\
               V(a)[X]^{j/2}, &\text{if $j$ even,}
  \end{cases}
\]
\item[V2:]\[V(b[X]^{k-1}d[X])=
   \begin{cases}(-1)^{|b|}\frac{1}{k}V((db)[X]^k)-(-1)^{|b|}\frac{2}{k}dV(b[X]^k), &\text{if $k$ odd,}\\
    V(b)[X]^{k/2-1}d[X], &\text{if $k$ even,} 
   \end{cases}\]
\item[V3:    ]$V(V^r(c[X]^l))=V^{r+1}(c[X]^l)),$
\item[V4:    ]$V(dV^s(e[X]^m))=2dV^{s+1}(e[X]^m)).$
\end{itemize}
The action of $F$ is given by: 
\begin{itemize}
\item[F1:]$F(a[X]^j)=F(a)[X]^{2j},$
\item[F2:]$F(b[X]^{k-1}d[X])=F(b)[X]^{2k-1}d[X],$
\item[F3:]$F(V^r(c[X]^l))=2V^{r-1}(c[X]^l),$
\item[F4:]$F(dV^s(e[X]^m))=dV^{s-1}(e[X]^m)+V^{s-1}(\iota(e)[X]^m.$
\end{itemize}
The action of $d$ is given by:
\begin{itemize}
\item[d1:]$d(a[X]^j)=d(a)[X]^j+(-1)^{|a|}ja[X]^{j-1}d[X],$
\item[d2:]$d(b[X]^{k-1}d[X])=d(b)[X]^{k-1}d[X]+k\iota(b)[X]^{k-1}d[X],$
\item[d3:]$d(V^r(c[X]^l)=dV^r(c[X]^l,$
\item[d4:]$d(dV^s(e[X]^m))=dV^s(\iota(e)[X]^m).$
\end{itemize}
On morphisms the functor $P$ is defined in the obvious way: if $\theta : E^*_{\bullet}\to F^*_{\bullet}$ is a morphisms of Witt complexes over $A$,  then $P(\theta):P(E^*_{\bullet})\to P(F^*_{\bullet})$ is defined on elements of type 1 by the formula $P(\theta)(a[X]^i)=\theta(a)[X]^i$ and similarly for elements of the other three types. 

\begin{thm}The functor $P:\mathcal{W}_A\to\mathcal{W}_{A[X]}$ is well defined and is a left adjoint of the forgetful functor $\mathcal{W}_{A[X]}\to\mathcal{W}_A$.
\end{thm}

\begin{proof}
The fact that the functor $P$ is well defined means that for any Witt complex $E$ over the ring $A$, the complex $P(E)$ is indeed a Witt complex over $A[X]$. We need to prove that the six conditions in the definition of a Witt complex are satisfied. 
Only two of these conditions and relations are hard to verify, the associativity and the relation $Fd\lambda([f]_n) = \lambda([f]_{n-1})d\lambda([f]_{n-1}]), \text{ for all } f\in A[X]$. The associativity requires a straightforward verification, that we will do in an appendix. 

We will prove the relation $Fd\lambda([f]_n) = \lambda([f]_{n-1})d\lambda([f]_{n-1}]), \text{ for all } f\in A[X]$ using induction by the level. For the level $n=1$ the identity is trivial, as both sides are equal to zero. Assume we know that the identity is true for the level $n-1$. We notice that the relation is easily verified for monomials, that is elements of the form $f=aX^m\in A[X]$.
\begin{lemma}The relation 
\[Fd\lambda([aX^m]_n) = \lambda([aX^m]_{n-1})d\lambda([aX^m]_{n-1})\]
holds for all $aX^m\in A[X]$.
\end{lemma}
\begin{proof}
Because there is no danger of confusion we will drop $\lambda$ and the subscript index indicating the level. 
\begin{align*}
Fd[aX^m]&=Fd([a][X]^m)=F((d[a])[X]^m+[a]d[X]^m)\\
&=F(d[a])F([X]^m)+F([a])F(d[X]^m)\\
&=([a]d[a])[X]^{2m}+m[a]^2F([X]^{m-1}d[X])\\
\intertext{and using the formula F2:}
Fd[aX^m]&=[a][X]^{2m}d[a]+m[a]^2[X]^{2m-1}d[X]\\
&=[a][X]^m([X]^md[a]+[a]d[X]^m)=[aX^m]d[aX^m],
\end{align*}
which is what we wanted to prove.
\end{proof}

The relation follows for arbitrary polynomials from the additivity result \ref{additivity}.

With this we proved that $P:\mathcal{W}_A\to\mathcal{W}_{A[X]}$ is well defined. To prove that it is the left adjoint of the forgetful functor $U:\mathcal{W}_{A[X]}\to\mathcal{W}_A$ we need to show that:
\[\text{Hom}_{\mathcal{W}_{A[X]}}(P(E),E')\cong\text{Hom}_{\mathcal{W}_{A}}(E, U(E')).\]
The morphism from left to right takes a map $f: P(E)\to E'$ to its restriction to $E\rightarrowtail U(P(E))$. The morphism from right to left takes $g :E\to U(E')$ to its unique extension $\tilde{g} :P(E)\to E'$ defined such that $g([X]_n)=\lambda'([X]_n).$ The two morphisms are inverse to each other. 
\end{proof}

\subsection{The 2-typical de Rham-Witt complex of the log-ring $\Z_{(2)}$ with the canonical log-structure}
In this section we define the notion of a 2-typical de Rham-Witt complex associated to a log-ring and we compute this complex for the ring $\Z_{(2)}$ with the canonical log-structure. We first recall the notions of log-rings and of differentials with log-structures. The standard reference is \cite{K}.
\begin{defn}
A log-ring is a ring $R$ together with a map of monoids
\[\alpha :M\to R, \]
where $R$ is considered a monoid under the multiplication. We will denote this log-ring be $(R,M)$. 
\end{defn}

The map $\alpha$ itself is called a ``pre-log structure''. 

\begin{defn}
A derivation of a log-ring $(R,M)$ into an $R$-module $E$ is a pair of maps
\[ (D, D\operatorname{log}):(R,M)\to E,\]
where $D: R\to E $ is a derivation and $D\operatorname{log}:M\to R$ is a map of monoids such that for all $a\in M$, $$\alpha(a)D\operatorname{log}(a)=Da.$$ 
\end{defn}

There is a universal example of a derivation of a log-ring $(R,M)$ given by the $R$-module 
\[\Omega^1_{(R,M)}=(\Omega^1_R\oplus(R\otimes_{\Z}M^{\operatorname{gp}}))/<d\alpha(a)-\alpha(a)\otimes a>,\]
where $M^{\operatorname{gp}}$ is the group completion of the monoid $M$. The structure maps are : 
\begin{align*}
d:R\to\Omega^1_{(R,M)},& \quad da=da\oplus0,\\
d\operatorname{log}:M\to\Omega^1_{(R,M)},& \quad d\operatorname{log}\, a =0\oplus(1\otimes a).
\end{align*}

\begin{defn}
A log-differential graded ring $(E^*,M)$ consists of a differential graded ring $E^*$, a pre-log structure $\alpha:M\to E^0$, and a derivation $(D,D\operatorname{log}):(E^0,M)\to E^1$ such that $D$ is equal to the differential $d:E^0\to E^1$ and such that $d\circ D\operatorname{log} =0$.
\end{defn}

The universal example of an anti-symmetric log-differential graded ring is:
\[\tilde{\Omega}^*_{(R,M)}=\tilde{\Lambda}_R^*(\Omega^1_{(R,M)}).\]
Here $\tilde{\Lambda}^*_R(N)=T_R(N)/\langle m\otimes n+n\otimes m \mid m, n\in N \rangle$ is the universal anti-symmetric graded $R$-algebra generated by the $R$-module $N$.

If $(R,M)$ is a log-ring, then for each $n\in \N$ the ring of length-$n$ Witt vectors, $W_n(R)$ over $R$ becomes part of the data that gives a log-ring $(W_n(R), M)$: the map of monoids $M \to W_n(R)$ is just the composition of the map $\alpha: M \to R$ and the Teichm\"{u}ller map $[-]_n:R\to W_n(R)$. 

\begin{defn}
A Witt complex $(E_{\bullet}^*, M_E)$ over a log-ring $(R,M)$ is a Witt complex $E_{\bullet}^*$ over $R$ together with pre-log structures $\,\alpha_n:M_E \to E_n^0$ and an extension of $\,\lambda: W_{\bullet}(R)\to E_{\bullet}^0$ to a strict map of pro-log-rings $\lambda: (W_{\bullet}(R),M)\to (E_{\bullet}^0,M_E),$ such that:
\begin{itemize}
\item[(i)]the pre-log structures $\,\alpha_n$ are compatible, in the sense that $R\circ \alpha_n=\alpha_{n-1}\,,$
\item[(ii)]$d\circ d\operatorname{log}[a]_n =0,$ for all $a\in M$,
\item[(iii)] $F\,d\operatorname{log}[a]_n=d\operatorname{log}[a]_{n-1},$ for all $a\in M.$
\end{itemize}
\end{defn}

\begin{prop}
The category of (2-typical) Witt complexes over a log-ring $(R,M)$ has an initial object $W_{\bullet}\Omega^*_{(R,M)}$, called the (2-typical) de Rham-Witt complex of $(R,M).$ 
\end{prop}
\begin{proof}
The proof is an application of the Freyd's adjoint functor theorem, entirely similar to the proof of Theorem \ref{initialobject} that asserts the existence of an initial object in the category of 2-typical Witt complexes over a ring.
\end{proof}
\textbf{Remark:} We note that if $M$ is the trivial monoid, then the 2-typical de Rham-Witt complex associated to $(R,M)$ is the 2-typical de Rham-Witt complex associated to $R,$ so the notion of a 2-typical de Rham-Witt complex is a generalization of the notion of a 2-typical de Rham-Witt complex. 
\\

In this section we will describe the 2-typical  de Rham-Witt complex of $(\Z_{(2)},M),$ where $M=\Q^*\cap\Z_{(2)}\hookrightarrow \Z_{(2)}$ is the canonical log-structure. The strategy is the same as in the previous calculations of de Rham-Witt complexes: we find a candidate $G_{\bullet}^*$ described explicitely by generators and relations and by formulas for the product and the actions of the various operators, and we prove that this candidate is isomorphic to $W_{\bullet}\Omega^*_{(R,M)}$. 

In degree zero the  de Rham-Witt complex is again the Witt vectors of the ring $R$, this following from a proof similar to the proof of Proposition \ref{degreezero}. In degree one, the only new generator that we have in the de Rham-Witt complex of $(\Z_{(2)},M),$ which is not in the de Rham-Witt complex of $\Z_{(2)}$ is $d\operatorname{log}[2]$. The product formulas are the same for the elements that already existed in the de Rham-Witt complex of $\Z_{(2)}$, so the only product formulas that we have to derive are $V^i(1)d\operatorname{log}[2]$.  

\begin{prop}i)The element $d\operatorname{log}[2]_n\in W_n\Omega^1_{(\Z_{(2)},M)}$ is annihilated by $2^n.$

ii) $V[1]_{n-1}d\operatorname{log}[2]_n\equiv 2d\operatorname{log}[2]_n\pmod{dV(W_{n-1}\Omega^0_{(\Z_{(2)},M)})}.$
\end{prop} 

\begin{proof}
The proof of both assertions is by induction. 
i)The case $n=1:$ $2d\operatorname{log}[2]_1=d(2)=0.$ Assuming $2^id\operatorname{log}[2]_i=0$ for all $i\leq n$ we will prove that $2^{n+1}d\operatorname{log}[2]_{n+1}=0.$ We use the formula $[2]_{n+1}d\operatorname{log}[2]_{n+1}=d[2]_{n+1}$ and the Corollary \ref{formulafor2} which says that 
\[[2]_{n+1}=\sum_{i=0}^nc_iV^i(1),\]
where $c_i=2^{-i}(2^{2^i}-2^{2^{i-1}}).$

We have:
\[\sum_{i=0}^nc_iV^i(1)d\operatorname{log}[2]_{n+1}=\sum_{i=1}^nc_idV^i(1)\]
We use that $V^i(1)d\operatorname{log}[2]_{n+1}=V^i(F^i(d\operatorname{log}[2]_{n+1}))=V^i(d\operatorname{log}[2]_{n+1-i}):$
\[2d\operatorname{log}[2]_{n+1}=\sum_{i=1}^n(-V^i(d\operatorname{log}[2]_{n+1-i})+c_idV^i(1)).\]
Now if we multiply this relation by $2^n$ we obtain: 
\[2^{n+1}d\operatorname{log}[2]_{n+1}=\sum_{i=1}^n(-V^i(2^{n+1}d\operatorname{log}[2]_{n+1-i})+2^{n+1}c_idV^i(1)),\]
which is zero by the induction hypothesis and the fact that $2^idV^i(1)=0.$

ii) For $n=1$ the congruence is trivial as both members are zero. Assuming that  the congruence holds for $n$, we want to prove that it holds for $n+1.$ Firs, by an easy induction, we see that:
\[V^i[1]_{n-i}d\operatorname{log}[2]_n\equiv 2^id\operatorname{log}[2]_n\pmod{dV}\]
Then we use the formula $[2]_{n+1}d\operatorname{log}[2]_{n+1}=d[2]_{n+1}$ combined with Corrolary \ref{formulafor2}: 
\[2d\operatorname{log}[2]_{n+1}+\sum_{i=1}^nc_iV^i[1]{n+1-i}d\operatorname{log}[2]_{n+1}=\sum_{i=1}^nc_idV^i[1]_{n+1-i}\equiv 0\pmod{dV}.\]
This gives:
\[2d\operatorname{log}[2]_{n+1}+\sum_{i=1}^nV(c_iV^{i-1}[1]{n-(i-1)}d\operatorname{log}[2]_{n})\equiv 0\pmod{dV},\]
or, using the formula described above for $V^{i-1}[1]{n-(i-1)}d\operatorname{log}[2]_{n}:$
\[2d\operatorname{log}[2]_{n+1}+V(\sum_{i=1}^nc_i2^{i-1}d\operatorname{log}[2]_{n})\equiv 0\pmod{dV},\]
\[2d\operatorname{log}[2]_{n+1}+V((2^{2^{n-1}}-1)d\operatorname{log}[2]_{n})\equiv 0\pmod{dV}.\]
Since $2^{2^{n-1}}$ is always divisible by $2^n,$ which annihilates $d\operatorname{log}[2]_{n},$ we get the desired result.

\end{proof}
   
The second part of this Proposition together with Proposition \ref{wittbasis}, which describes a basis of  $W_n(\field{Z}),$ tell us  that the product formula we are trying to derive is of the form:        
\[V[1]_nd\operatorname{log}[2]_n= 2d\operatorname{log}[2]_n+a_1dV[1]_{n-1}+a_2dV^2[1]_{n-2}+\cdots+
a_{n-1}dV^{n-1}[1]_1.\]

We note that the coefficients $a_i,$ don't depend on $n$, as we can apply $R$ to the relation in level $n+1$ to obtain the relation in level $n.$

\begin{lemma}
Assuming the previous product formula, the following formulas hold: 
\[V^i(1)d\operatorname{log}_n(2)=2^id\operatorname{log}_n(2)+2^{i-1}(a_1+\cdots+a_i)dV^i(1)+\cdots+2^{i-1}(a_{n-i}+\cdots+a_{n-1})dV^{n-1}(1).\]
\end{lemma}

\begin{proof}
The proof is by induction. The case $i=1$ is obvious. We prove that if the formula is true for $i$ then it must be true for $i+1$. 

\begin{align*}
V^{i+1}(1)d\operatorname{log}_n(2)&=V(V^i(1)d\operatorname{log}_{n-1}(2))\\
&=V(2^id\operatorname{log}_{n-1}(2)+2^{i-1}(a_1+\cdots+a_i)dV^i(1)+\cdots\\
&\phantom{=}\cdots +2^{i-1}(a_{n-i-1}+\cdots+a_{n-2})dV^{n-2}(1))\\
&=2^iV(1)d\operatorname{log}_n(2)+2^{i-1}(a_1+\cdots+a_i)VdV^i(1)+\\
&\phantom{=}\cdots+2^{i-1}(a_{n-i-1}+\cdots+a_{n-2})VdV^{n-2}(1)\\
&\hspace{-1.5cm}=2^{i+1}d\operatorname{log}_n(2)+2^ia_1dV(1)+\cdots+2^ia_idV^i(1)+2^ia_{i+1}dV^{i+1}(1)+\cdots\\
&\hspace{-1.5cm}\phantom{=}\cdots+2^ia_{n-1}dV^{n-1}(1)+2^{i}(a_1+\cdots+a_i)dV^{i+1}(1)+\cdots\\
&\hspace{-1.5cm}\phantom{=}\cdots+2^{i}(a_{n-i-1}+\cdots +a_{n-2})dV^{n-1}(1)\\
&\hspace{-2.2cm}=2^{i+1}d\operatorname{log}_n(2)2^{i}(a_1+\cdots+a_{i+1})dV^{i+1}(1)+\cdots+2^{i}(a_{n-i-1}+\cdots+a_{n-1})dV^{n-1}(1).
\end{align*}

We used the fact that $2^idV(1)=\cdots=2^idV^i(1)=0$, which follows from $2dV=Vd$ and $d(1)=0$.
\end{proof}

We will now compute the coefficients $a_i$. We start with the relation
\[[2]_nd\operatorname{log}[2]_n=d[2]_n\]
This gives: 
\[\sum_{i=0}^{n-1}c_i(V^i(1)d\operatorname{log}_n(2))=\sum_{i=1}^{n-1}c_idV^i(1),\]
and we use the formula that we just derived for $V^i(1)d\operatorname{log}_n(2)$: 
\[c_0d\operatorname{log}_n(2)+\sum_{i=0}^{n-1}c_i(2^id\operatorname{log}_n(2)+\sum_{j=i}^{n-1}(\sum_{k=j-i+1}^ja_k)dV^j(1))=\sum_{i=1}^{n-1}c_idV^i(1).\]

We regroup the sums and we obtain: 
\[(\sum_{i=0}^{n-1}2^ic_i)d\operatorname{log}_n(2)+\sum_{j=1}^{n-1}\sum_{i=1}^j(2^ic_i\sum_{k=j-i+1}^ja_k)dV^j(1)=\sum_{j=1}^{n-1}c_jdV^j(1).\]
The first term in the left hand side member is zero because $\sum_{i=0}^{n-1}2^ic_i=2^{2^{n-1}}$ and $d\operatorname{log}_n(2)$ is annihilated by $2^n$. We equate the coefficients of $dV^j(1)$ modulo $2^j$ and obtain:
\[\sum_{i=1}^j(2^ic_i\sum_{k=j-i+1}^ja_k)\equiv c_j\pmod{2^j},\]
or: 
\[\sum_{k=1}^j(\sum_{i=j-k+1}^j2^{i-1}c_i)a_k\equiv c_j\pmod{2^j}.\]
Let us call $B_{jk}=\sum_{i=j-k+1}^j2^{i-1}c_i,$ if $j\geq k$, and $B_{jk}=0,$ if$j<k$. We obtain therefore a system of equations: 
$$\sum_{k=1}^jB_{jk}a_k\equiv c_j \pmod{2^j}.$$ 
We need to make a comment about this system. A priori the unknowns $a_k$ are in different rings, namely $a_k\in \Z/{2^k\Z}$. So the system as it stands doesn't really make sense. However we can think of $2^{j-k}a_k$ as an element of $\Z/{2^j\Z}$, and we notice that $B_{jk}$ is divisible by $2^{j-k}$, because $B_{jk}=\sum_{i=j-k+1}^j2^{i-1}c_i=2^{2^{j}-1}-2^{2^{j-k}-1},$ if $j\geq k$, and $B_{jk}=0,$ if$j<k$.  To solve this system, we lift it to a system over the ring of integers $\Z$, we solve that system, and then take classes of congruence modulo the corresponding power of 2. 

We observe that $B_{jk}=2^{2^{j}-1}-2^{2^{j-k}-1}\equiv -2^{2^{j-k}-1} \pmod{2^j}$, and that $c_1=1$, $c_2\equiv -1\pmod{2^2}$, $c_3\equiv-2\pmod{2^3}$, and $c_i\equiv 0 \pmod{2^i},$ for $i>3$. We find thus convenient to lift the previous system of equation to the following system over $\Z$:
$$\sum_{k=1}^jb_{jk}a_k=c_j' ,$$
with $b_{jk}=-2^{2^{j-k}-1}$ for $j\geq k$, $b_{jk}=0$ for $j< k$, $c_1'=1$, $c_2'= -1$, $c_3'=-2$, and $c_i= 0,$ for $i>3$.   

The matrix of the system is lower triangular and it has only $-1$ on the diagonal. We can invert it using, for example, Gauss-Jordan's method. The inverse matrix $F$ is also lower triangular and it has entries: 
\[f_{ij}=
\begin{cases}
0,&\text{if $i<j,$}\\
-1,&\text{if $i=j,$}\\
\sum_{i=i_0>i_1>\cdots>i_s=j}b_{i_0i_1}b_{i_1i_2}\cdots b_{i_{s-1}i_s}, &\text{if $i> j.$}
\end{cases} \]
A direct computation shows that $b_{21}=b_{32}=b_{43}=-2,$ $b_{31}=b_{42}=-2^3,$ and using this that $f_{21}=f_{32}=f_{43}=2$ and $f_{31}=f_{42}=-2^2.$ These entries of the matrix $F$ are all we need to compute the first three coeficients in the product formula. We obtain: $a_1\equiv 1\pmod2, a_2\equiv -1\pmod4, a_3\equiv 4 \pmod8$. We will prove that all the other coefficients are zero. 

For all $i\geq 3$ we have: 
\[a_i\equiv f_{i1}c_1'+ f_{i2}c_2'+f_{i3}c_3'= f_{i1}-f_{i2}-2f_{i3}\pmod{2^i}\]     
We remark first that $v_2(b_{ij})\geq i-j$, with equality if and only if $i=j+1.$ Using this observation we see that all the terms that add up to give $f_{i1}$ are  divisible by $2^{i-1}$, and the only one that is not divisible by $2^i$ is $b_{ii-1}b_{i-1i-2}\cdots b_{21}=(-1)^{i-1}2^{i-1}$. Therefore $f_{i1}\equiv (-1)^{i-1}2^{i-1}\pmod{2^i}.$ Similarly, $f_{i2}$ is divisible by $2^{i-2}$ and the terms in the sum that makes up $f_{i2}$ that are not divisible by $2^i$ are of two forms: 
\begin{itemize}
\item[-]one product $b_{ii-1}b_{i-1i-2}\cdots b_{32}=(-1)^{i-2}2^{i-2},$\\
\item[-]$(i-3)$ products of the form $b_{ii-1}b_{i-1i-2}\cdots b_{k+1k}b_{kk-2}b_{k-2k-3}\cdots b_{32}=(-1)^{i-3}2^{i-1}.$
\end{itemize}
We obtain $f_{i2}\equiv (-1)^{i-2}2^{i-2}+(i-3)(-1)^{i-3}2^{i-1}\pmod{2^i}.$

We treat $f_{i3}$ in the same way: 
\begin{align*}2f_{i3}&\equiv 2(b_{ii-1}b_{i-1i-2}\cdots b_{32}+\sum_kb_{ii-1}b_{i-1i-2}\cdots b_{k+1k}b_{kk-2}b_{k-2k-3}\cdots b_{32}\\
&=2[(-1)^{i-3}2^{i-3}+(i-4)2^{i-2}]\pmod{2^i}
\end{align*}
With these formulas we can compute $a_i$ for $i>3$:
\begin{align*}
a_i&=f_{i1}-f_{i2}-2f_{i3}\\
&\equiv (-1)^{i-1}2^{i-1}-(-1)^{i-2}2^{i-2}-(i-3)(-1)^{i-3}2^{i-1}-2[(-1)^{i-3}2^{i-3}+(i-4)2^{i-2}]\\
&\equiv 0 \pmod{2^i}
\end{align*}
We have:
\begin{lemma}
The following product formula holds for all $n$: 
\[V[1]_nd\operatorname{log}[2]_n= 2d\operatorname{log}[2]_n+dV[1]_{n-1}-dV^2[1]_{n-2}+4dV^3[1]_{n-3}.\]
\end{lemma}

Now we can state the structure theorem for the 2-typical de Rham-Witt complex of $(\Z_{(2)},M),$ where $M=\Z_{(2)}^*.$

\begin{thm}
The structure of  $W_{\bullet}\Omega^*_{(\Z_{(2)},M)}$ is:

\begin{itemize}
\item[(i)]As abelian groups
\begin{align}
W_n\Omega^0_{(\Z_{(2)},M)}&=\bigoplus_{i=0}^{n-1}\Z_{(2)}\cdot V^i(1),\\
W_n\Omega^1_{(\Z_{(2)},M)}&=\Z/{2^n\Z}\,d\operatorname{log}[2]_n \oplus\bigoplus_{i=1}^{n-1}\Z/{2^i\Z}\cdot dV^i(1),\\
W_n\Omega^i_{(\Z_{(2)},M)}&=0, \text{for all } i\geq2.
\end{align}

The product relations and the actions of the various operators are the ones from Theorem \ref{structureforZ}, in addition to which we also have:

\item[(ii)]the product formulas: 
\[V[1]_{n-1}d\operatorname{log}[2]_n= 2d\operatorname{log}[2]_n+dV[1]_{n-1}-dV^2[1]_{n-2}+4dV^3[1]_{n-3},\]
\[V^i[1]_{n-i}d\operatorname{log}[2]_n= 2^id\operatorname{log}[2]_n-2^{i-1}dV^{i+1}[1]_{n-i-1}+2^{i+1}dV^{i+2}[1]_{n-i-2},\]
\item[(iii)]the action of the operator $V$:
\[V(d\operatorname{log}[2]_n)= 2d\operatorname{log}[2]_{n+1}+dV[1]_{n}-dV^2[1]_{n-1}+4dV^3[1]_{n-2},\]
\item[(iv)]the action of the operator $F$:
\[F(d\operatorname{log}[2]_n)=d\operatorname{log}[2]_{n-1}.\]
\end{itemize}
\end{thm}
\begin{proof}
The proof is similar to the proof of Theorem \ref{structureforZ}, which describes the structure of the 2-typical de Rham-Witt complex of the integers. 
More precisely we define the pro-graded ring $G_{\bullet}^*$ to be:
 \begin{align}
G_n^0&= \bigoplus_{i=0}^{n-1}\Z_{(2)}\cdot V^i(1),\\
G_n^1&= \Z/{2^n\Z}\,d\operatorname{log}[2]_n\oplus\bigoplus_{i=1}^{n-1}\Z/{2^i\Z}\cdot dV^i(1),\\
G_n^i&=0,  \text{for all } i\geq2,
\end{align}
with the product rule and the action of the operators as in the theorem. We want to prove that $W_{\bullet}\Omega^*_{(\Z_{(2)},M)}\cong G_{\bullet}^*.$ We will show that we have morphisms
\begin{align*}
\phi&:W_{\bullet}\Omega^*_{(\Z_{(2)},M)}\to G_{\bullet}^*,\\
\psi&:G_{\bullet}^*\to W_{\bullet}\Omega^*_{(\Z_{(2)},M)},
\end{align*} 
and that they are inverse to each other. 

The existence (and uniqueness) of the morphism $\phi$ follows from the fact that $G_{\bullet}^*$ is a 2-typical Witt complex. The definition of $\psi$ is forced by the requirements that it is a morphism of 2-typical Witt complexes: $[1]_n\mapsto [1]_n$ (since $\psi$ is a morphism of rings), $V^i[1]_{n-i}\mapsto V^i[1]_{n-i},$ $dV^i[1]_{n-i}\mapsto dV^i[1]_{n-i}$ ($\psi$ commutes with $V$ and $d$), $[2]_n\mapsto [2]_n$ (because $[2]_n=\sum_{i=0}^{n-1}c_iV^i[1]_{n-i}$ by Prop \ref{formulafor2} and $\psi$ is additive), $d\operatorname{log}[2]_n\mapsto d\operatorname{log}[2]_n$ (because $\psi$ commutes with $d\operatorname{log}$).

In order to see that $\psi$ is well defined we need to show that $W_n\Omega^i_{(\Z_{(2)},M)}=0$ for $i\geq2.$ This is proven by induction on $n$. The first step of the iduction, that $W_1\Omega^i_{(\Z_{(2)},M)}=0$ follows from the fact that $\Omega^i_{(\Z_{(2)},M)}=0$ for all $i\geq 2$ (which follows from $d\circ d\operatorname{log} =0$) and the fact that $\Omega^i_{(\Z_{(2)},M)}\to W_1\Omega^i_{(\Z_{(2)},M)}$ is surjective. 

Assuming that $W_{n-1}\Omega^i_{(\Z_{(2)},M)}=0$ for $i\geq2,$ we want to prove that $W_n\Omega^i_{(\Z_{(2)},M)}=0$ for $i\geq2.$ We use the standard filtration:
\[\operatorname{Fil}^sW_n\Omega^i_{(\Z_{(2)},M)}=V^sW_{n-i}\Omega^i_{(\Z_{(2)},M)}+dV^sW_n\Omega^{i-1}_{(\Z_{(2)},M)}. \]
The sequence:
\[0\to\operatorname{Fil}^{n-1}W_n\Omega^i_{(\Z_{(2)},M)}\to W_n\Omega^i_{(\Z_{(2)},M)}\to W_{n-1}\Omega^i_{(\Z_{(2)},M)}\to 0\] 
is exact by exactly the same argument used to prove Lemma \ref{exactsequence}. For $i\geq 2$ the last term in this short exact sequence is zero by the induction hypothesis. The first term is: 
\begin{align*}\operatorname{Fil}^{n-1}W_n\Omega^i_{(\Z_{(2)},M)}&=V^{n-1}W_{n-i}\Omega^i_{(\Z_{(2)},M)}+dV^{n-1}W_n\Omega^{i-1}_{(\Z_{(2)},M)}\\
&=V^{n-1}(0)+dV^{n-1}(\Z/2\Z\,d\operatorname{log}[2]_1).
\end{align*}
This is zero if $dV^{n-1}(d\operatorname{log}[2]_1)=0.$ We have:
\begin{align*}
dV^{n-1}(d\operatorname{log}[2]_1)&=d(V^{n-1}([1]_1)d\operatorname{log}[2]_n)\\
&=d(2^{n-1}d\operatorname{log}[2]_n+dV(x))\\
&=2^{n-1}d\circ d\operatorname{log}[2]_n+ddV(x)\\
&=ddV(x),
\end{align*}
where $x\in W_{n-1}\Omega^{i-2}_{(\Z_{(2)},M)}.$ Too see that $ddV(x)=0$ we use the following trick: $W_{\bullet}\Omega^*_{(\Z_{(2)},M)}$ is a Witt complex over the log-ring $(\Z_{(2)},M),$ so in particular it is a Witt complex over the ring $\Z_{(2)},$ and as such it is the target of a unique homomorphism from the de Rham-Witt complex $W_{\bullet}\Omega^*_{\Z_{(2)}}.$
The element $ddV(x)$ is in the image of this homomorphism, but $W_{\bullet}\Omega^i_{\Z_{(2)}}=0$ for $i\geq 2$, therefore  $ddV(x)=0.$ This finishes the proof that $W_{n}\Omega^{i}_{(\Z_{(2)},M)}=0$ if $i\geq 2$, and thus $\psi$ is well defined. 

Too see that $\phi$ and $\psi$ are inverse to each other we check that $\phi\circ\psi=1$ and $\psi\circ\phi=1.$ The first follows from the fact that $\phi\circ\psi$ is a morphism of Witt complexes, and therefore $[1]_n\mapsto [1]_n,$ $V^i[1]_{n-i}\mapsto V^i[1]_{n-i},$ $dV^i[1]_{n-i}\mapsto dV^i[1]_{n-i},$ $[2]_n\mapsto [2]_n,$ and $d\operatorname{log}[2]_n\mapsto d\operatorname{log}[2]_n,$ and the second from the fact that $\psi\circ\phi$ is an endomorphism of an initial object in a category.
\end{proof}

\section{Appendix A: The $p$-typical de Rham-Witt complex of $\Z_{(p)}$ with the canonical log-structure, $p$ odd}

In this appendix we give the structure of the $p$-typical de Rham-Witt complex of the log-ring $(\Z_{(p)},M),$ for $p$ odd. Here $M=\Q^*\cap \Z_{(p)}\hookrightarrow  \Z_{(p)}$ is the canonical log-structure on  $\Z_{(p)}.$  The computations are exactly the same as for the case $p=2$, but the results are a little different. We first recall the structure of the $p$-typical de Rahm-Witt complex of $\Z_{(p)}$ from Example 1.2.4 of \cite{He2}. 

\begin{prop} The structure of $W_{\bullet}\Omega^*_{\Z_{(p)}}$ for $p$ odd is:
\begin{itemize}
\item[(i)]As abelian groups: 
\begin{eqnarray}
W_n\Omega^0_{\Z_{(p)}}&=& \bigoplus_{i=0}^{n-1}\Z_{(p)}\cdot V^i(1),\\
W_n\Omega^1_{\Z_{(p)}}&=& \bigoplus_{i=1}^{n-1}\Z/{p^i\Z}\cdot dV^i(1) ,\\
W_n\Omega^i_{\Z{(p)}}&=&0,\quad \mbox{for } i\geq 2.
\end{eqnarray}
\item[(ii)] The product is given by:
\begin{align}
V^i(1)V^j(1)&=p^iV^j(1), \text{if } i\leq j,\\
V^i(1)dV^j(1)&=\begin{cases}
                 p^idV^j(1),& \text{if } i<j,\\
                 0, &\text{if } i\geq j.
               \end{cases}
\end{align}
\item[(ii)] The action of the operators $F$ and $V$ is given by:
\begin{align}
FV^i(1)&=pV^{i-1}(1),\\
FdV^i(1)&=dV^{i-1}(1),\\
V(V^i(1))&=V^{i+1}(1),\\
V(dV^i(1))&=pdV^{i+1}(1).
\end{align}
\end{itemize}
\end{prop}

The structure of $W_\bullet\Omega^*_{(\Z_{(p)},M)}$ is different for $p=3$ and $p\geq 5$. 

\begin{thm} The structure of $W_\bullet\Omega^*_{(\Z_{(p)},M)}$ with $p$ odd is:
\begin{itemize}
\item[(i)]As additive groups: 
\begin{align}
W_n\Omega^0_{(\Z_{(p)},M)}&= W_n\Omega^0_{\Z_{(p)}}\cong W_n({\Z_{(p)}}),\\
 W_n\Omega^1_{(\Z_{(p)},M)}&= W_n\Omega^1_{\Z_{(p)}}\oplus \Z/{p^n\Z}d\operatorname{log}[p]_n,\\
W_n\Omega^i_{(\Z_{(p)},M)}&=0, \text{ for all } i\geq 2.
\end{align}
\item[(ii)] The product is given by the formulas in the previous theorem and the following formula that involves $d\operatorname{log}[p]$:
\[V^i[1]_{n-i}d\operatorname{log}[p]_n=
\begin{cases}
3^id\operatorname{log}[3]_n+3^{i-1}dV^i[1]_{n-i}+3^idV^{i+1}[1]_{n-i-1},& \text{if } p=3,\\
p^i\operatorname{log}[p]_n+p^{i-1}dV^i[1]_{n-i},& \text{if } p\geq 5.
\end{cases}\]
\end{itemize}
\item[(ii)] The action of $F$ and $V$ on $d\operatorname{log}[p]_n$ is:
\begin{align}
F(d\operatorname{log}[p]_n)&=d\operatorname{log}[p]_{n-1},\\
V(d\operatorname{log}[p]_n)&=
\begin{cases}
3d\operatorname{log}[3]_{n+1}+dV[1]_n+3dV[1]_{n-1},& \text{if } p=3,\\
p\,d\operatorname{log}[p]_{n+1}+dV[1]_n,& \text{if } p\geq 5.
\end{cases}
\end{align}
\end{thm}
The proof of this theorem is entirely similar to the proof of the structure theorem for $W_\bullet\Omega^*_{(Z_{(2)},M)}.$ 

\section{Appendix B: Associativity }

In this appendix we discuss the associativity of the multiplication rule defined in Section 4. We recall that the functor $P:\mathcal{W}_A\to\mathcal{W}_{A[X]}$ is defined on objects as follows: for a a Witt complex $E^*_{\bullet}\in \mathcal{W}_A$, $P(E)_n^q$ consists of formal sums of four types of elements: 
 \begin{itemize}
\item Type 1: elements of the form $a[X]^j$, where $a\in E_n^q$,
\item Type 2: elements of the form $b[X]^{k-1}d[X]$, where $b\in E_n^{q-1}$,
\item Type 3: elements of the form $V^r(c[X]^l)$, where $r>0$, $c\in E_{n-r}^q$, and $l$ is odd,
\item Type 4: elements of the form $dV^s(e[X]^m)$, where $s>0$, $e\in E_{n-s}^{q-1}$, and $m$ is odd.
\end{itemize}
The product is given by ten formulas, from P1.1 to P1.4. 

We make now the convention that, for example, A1.3.4 means the statement that says that $(xy)z=x(yz),$ where $x$ is an element of the first type, $y$ an element of the third type, and $z$ an element of the fourth type. In order to prove the associativity one has to check twenty relations like this, from A1.1.1 to A4.4.4\,. 

Since there are three product formulas given in ``cases'' format, the associativity relations involving these formulas will be a little more tedious to verify. Ten out of the twenty relations that we want to check contain at least a product given in cases format. Out of the remaining ten, five are more or less trivial, namely A1.1.1, A1.1.2, A1.1.3, A1.2.2, and A2.2.2\,. The five formulas that don't involve products with the cases format are A1.1.4, A1.2.3, A1.2.4, A2.2.3, and A2.2.4\,. The hardest seems to be the first, even if it doesn't involve the operator $\iota.$ We will show how it is derived, and then we will also show A1.2.4\,, where $\iota$ is involved. 

Among the ten cases where at least one product is in the cases format, one is almost trivial, A1.3.3\,. The other nine require all some extensive computations. The most dificult of them are A3.3.4 and A3.4.4\,.We will show how A3.3.4 is derived. 

We start with the relation A1.1.4\,. Let $x=aX^j,$ $y=a'X^{j'},$ and $z=dV^s(eX^m).$ Then:
\begin{align*}
(xy)z&=aa'X^{j+j'}dV^s(eX^m)\\
&=(-1)^{|aa'|}\frac{m}{2^s(j+j')+m}dV^s(F^s(aa')eX^{2^s(j+j')+m})\\
&\hspace{-1cm}-(-1)^{|aa'|}V^s((F^s(d(aa'))e-\frac{j+j'}{2^s(j+j')+m}d(F^s(aa')e))X^{2^s(j+j')+m}).
\end{align*}
On the other hand:
\begin{align*}
x(yz)=&aX^j\{(-1)^{|a'|}\frac{m}{2^sj'+m}dV^s(F^s(a')eX^{2^sj'+m})\\
&-(-1)^{|a'|}V^s((F^s(da')e-\frac{j'}{2^sj'+m}d(F^s(a')e))X^{2^sj'+m})\}.
\end{align*}
If we denote $E=F^s(a')e,$ $M=2^sj'+m,$ and $P=2^{s}(j+j')+m,$ we obtain:
\begin{align*}
x(yz)=&aX^j\{(-1)^{|a'|}\frac{m}{M}dV^s(EX^M)-(-1)^{|a'|}V^s((E-\frac{j'}{M}d(F^s(a')e))X^M)\}\\
=&(-1)^{|a'|}\frac{m}{M}\{(-1)^{|a|}\frac{M}{P}dV^s(F^s(a)EX^P)\\
&-(-1)^{|a|}V^s((F^s(da)E-\frac{j}{P}d(F^s(a)E))X^P)\}\\
&-(-1)^{|a'|}V^s((F^s(ada')e-\frac{j'}{M}F^s(a)d(F^s(a')e))X^P)\\
=&(-1)^{|aa'|}\frac{m}{P}dV^s(F^s(aa')eX^{P})+V^s(UX^P),
\end{align*}
where $U$ is the expression:
\begin{align*}
U=&(-1)^{|aa'|}\frac{m}{M}F^s(da\,a')e+(-1)^{|aa'|}\frac{m}{M}\frac{j}{P}d(F^s(aa')e)\\
&-(-1)^{|a'|}F^s(ada')e+(-1)^{|a'|}F^s(a)d(F^s(a')e)).
\end{align*}
Using the fact that $d$ is a derivation and that $dF^s=2^sF^sd,$ we obtain:
\begin{align*}
U=&-(-1)^{|aa'|}\frac{m}{P}F^s(da\,a')e-(-1)^{|a'|}\frac{m}{P}F^s(ada')e\\
&+\frac{j+j'}{P}F^s(aa')de,
\end{align*}
which agrees with the expression we find inside $V^s$ for the elemtent $(xy)z.$

We prove now A1.2.4\,. Let $x=aX^j,$ $y=bX^{k-1}dX,$ and $z=dV^s(eX^m).$ We have: 
\begin{align*}
(xy)z=&(abX^{j+k-1}dX)dV^s(eX^m)\\
=&-(-1)^{|ab|}\frac{1}{2^s(j+k)+m}V^s(F^s(d(ab)+(j+k)\iota(ab)deX^{2^s(j+k)+m}))\\
&+(-1)^{|ab|}\frac{1}{2^s(j+k)+m}dV^s(F^s(ab)deX^{2^s(j+k)+m}).
\end{align*}
On the other hand: 
\begin{align*}
x(yz)=&aX^j\{-(-1)^{|b|}\frac{1}{2^sk+m}V^s(F^s(db+(j+k)\iota(b)deX^{2^sk+m}))\\
&+(-1)^{|b|}\frac{1}{2^sk+m}dV^s(F^s(b)deX^{2^sk+m})\}.
\end{align*}
We denote by $M=2^sk+m,$ and $P=2^s(j+k)+m,$ we obtain:
\begin{align*}
x(yz)=&-(-1)^{|b|}\frac{1}{M}V^s(F^s(a\,db+k\iota(ab))deX^P)\\
&+(-1)^{|b|}\frac{1}{M}\{(-1)^{|a|}\frac{M}{P}dV^s(F^s(ab)deX^P)\\
&-(-1)^{|a|}V^s((F^s(da)F^s(b)de-\frac{j}{P}d(F^s(ab)de))X^P)\}\\
=&(-1)^{|ab|}\frac{1}{P}dV^s(F^s(ab)deX^P)+V^s(WX^P).
\end{align*}
We compute separately the expression $W:$
\begin{align*}
W=&-(-1)^{|b|}\frac{1}{M}F^s(a\,db)de+kF^s(\iota(ab))de\\
&-(-1)^{|ab|}\frac{1}{M}F^s(da\,b)de+(-1)^{|ab|}\frac{j}{MP}d(F^s(ab)de)\\
=&-(-1)^{|ab|}\frac{1}{M}F^s(d(ab))de+F^s(k\iota(ab))de\\
&+(-1)^{|ab|}\frac{2^sj}{MP}F^s(d(ab))de+\frac{j}{MP}F^s(ab)dde.
\end{align*}
The last term in this sum is $\frac{j}{MP}F^s(ab)dde=\frac{j}{MP}F^s(\iota(ab))de.$ We make the observation that for any odd number $m\in \Z$ we have $\frac{1}{m}\iota=\iota,$ so the last term in the sum becomes simply $jF^s(\iota(ab))de.$ Therefore:
\begin{align*}
W=&-(-1)^{|ab|}\frac{1}{P}F^s(d(ab))de+F^s((j+k)\iota(ab))de\\
=&-(-1)^{|ab|}\frac{1}{P}F^s(d(ab)+(j+k)\iota(ab))de,
\end{align*}
which agrees with the expression we find inside $V^s$ in the product $(xy)z.$

We show now how to derive the associativity relation A3.3.4\,. Let $x=V^r(cX^l),$ $y=V^{r'}(c'X^{l'}),$ and $z=dV^s(eX^m).$ Since the product formulas depend on the ordering of the exponents $r,$ $r',$ and $s,$ it follows that verifying this relation involves checking 13 different cases, from $r<r'<s$ to $s<r'<r.$ We will verify the case $s<r'<r.$ 

We make the notations  $C=cF^{r-r'}(c'),$ $C'=c'F^{r'-s}((d+\iota)(e)),$ $E=c'F^{r'-s}(e),$ $L=2^{r-r'}l'+l,$ $M=2^{r'-s}m+l',$ and $P=2^{r-s}m+2^{r-r'}l'+l.$ We have:
\begin{align*}
(xy)z=&(V^r(cX^l)V^{r'}(c'X^{l'}))dV^s(eX^m)\\
=&2^{r'}V^r(cF^{r-r'}(c')X^{2^{r-r'}l'+l})dV^s(eX^m)\\
=&2^{r'}V^r(CX^L)dV^s(eX^m)\\
=&2^{r'}\{V^r(CF^{r-s}((d+\iota)(e))X^P)+(-1)^{|cc'|}\frac{2^rm}{P}dV^r(CF^{r-s}(e)X^P)\\
&-(-1)^{|cc'|}\frac{m}{P}V^r(d(CF^{r-s}(e))X^P)\}\\
=&(-1)^{|cc'|}\frac{2^{r+r'}m}{P}dV^r(cF^{r-r'}(c')F^{r-s}(e)X^P)+V(UX^P),
\end{align*}
where the expression $U$ is:
\begin{align*}
U=&2^{r'}\{cF^{r-r'}(c')F^{r-s}((d+\iota)(e))-(-1)^{|cc'|}\frac{m}{P}d(cF^{r-r'}(c')F^{r-s}(e))\}
\end{align*}

On the other hand: 
\begin{align*}
x(yz)=&V^r(cX^l)(V^{r'}(c'X^{l'})dV^s(eX^m))\\
=&V^r(cX^l)\{V^{r'}(c'F^{r'-s}((d+\iota)(e))X^{2^{r'-s}+l'})\\
&+(-1)^{|c'|}\frac{2^{r'}m}{2^{r'-s}m+l'}dV^{r'}(c'F^{r'-s}(e)X^{2^{r'-s}+l'})\\
&-(-1)^{|c'|}\frac{m}{2^{r'-s}m+l'}V^{r'}(d(c'F^{r'-s}(e))X^{2^{r'-s}+l'})\}\\
=&V^r(cX^l)\{V^{r'}(C'X^M)+(-1)^{|c'|}\frac{2^{r'}m}{M}dV^{r'}(EX^M)\\
&-(-1)^{|c'|}\frac{m}{M}V^{r'}(dE\,X^M)\}\\
=&2^{r'}V^r(cF^{r-r'}(C')X^P)+(-1)^{|c'|}\frac{2^{r'}m}{M}\{V^r(cF^{r-r'}((d+\iota)E)X^P)\\
&+(-1)^{|c|}\frac{2^{r}M}{P}dV^r(cF^{r-r'}(E)X^P)-(-1)^{|c|}\frac{M}{P}V^r(d(cF^{r-r'}(E))X^P)\}\\
&-(-1)^{|c'|}\frac{m}{M}2^{r'}V^r(cF^{r-r'}(dE)X^P).
\end{align*}
The second and the fifth term in this sum cancel each other, since $2^{r'}\iota=0.$ We have:
\[x(yz)=(-1)^{|cc'|}\frac{2^{r+r'}m}{P}dV^r(cF^{r-r'}(c')F^{r-s}(e)X^P)+V(WX^P),\]
where the expression $W$ is:
\begin{align*}
W=&2^{r'}\{cF^{r-r'}(c')F^{r-s}(d(e))+(-1)^{|c'|}\frac{2^{r'}m}{M}cF^{r-r'}(d(c'F^{r'-s}(e)))\\
&-(-1)^{|cc'|}\frac{m}{P}d(cF^{r-r'}(c'F^{r'-s}(e)))-(-1)^{|c'|}\frac{2^{r'}m}{M}cF^{r-r'}(d(c'F^{r'-s}(e))).
\end{align*}
The second and the fourth term cancel, and we see that $W=U$, hence $(xy)z=x(yz).$

\vspace{5mm}
\noindent {\small \sc 280 Marin Blvd,\\
Jersey City, NJ 07302, USA} \\

\noindent {{\em E-mail address}: Viorel.Costeanu@gmail.com}

\end{document}